
\def\LaTeX{%
  \let\Begin\begin
  \let\End\end
  \let\salta\relax
  \let\finqui\relax
  \let\futuro\relax}

\def\UK{\def\our{our}\let\sz s}
\def\USA{\def\our{or}\let\sz z}

\UK



\LaTeX

\USA


\salta

\documentclass[twoside,12pt]{article}
\setlength{\textheight}{24cm}
\setlength{\textwidth}{16cm}
\setlength{\oddsidemargin}{2mm}
\setlength{\evensidemargin}{2mm}
\setlength{\topmargin}{-15mm}
\parskip2mm



\usepackage[usenames,dvipsnames]{color}
\usepackage{amsmath}
\usepackage{amsthm}
\usepackage{amssymb}
\usepackage[mathcal]{euscript}
\usepackage{enumitem}
%
%
\usepackage{cite}
%
%
%


\definecolor{viola}{rgb}{0.3,0,0.7}
\definecolor{ciclamino}{rgb}{0.5,0,0.5}
\definecolor{rosso}{rgb}{0.8,0,0}

\def\andrea #1{{\color{red}#1}}
\def\dafare #1{{\color{red}#1}}
\def\andrea #1{#1}
\def\dafare  #1{#1}


\bibliographystyle{plain}


%

\finqui

\def\Beq{\Begin{equation}}
\def\Eeq{\End{equation}}
\def\Bsist{\Begin{eqnarray}}
\def\Esist{\End{eqnarray}}

\def\Bthm{\Begin{theorem}}
\def\Ethm{\End{theorem}}
\def\Blem{\Begin{lemma}}
\def\Elem{\End{lemma}}

\def\Bcor{\Begin{corollary}}
\def\Ecor{\End{corollary}}
\def\Brem{\Begin{remark}\rm}
\def\Erem{\End{remark}}

\def\Bcenter{\Begin{center}}
\def\Ecenter{\End{center}}
\let\non\nonumber




\def\step #1 \par{\medskip\noindent{\bf #1.}\quad}


\def\Lip{Lip\-schitz}
\def\Holder{H\"older}
\def\Frechet{Fr\'echet}

\def\aand{\quad\hbox{and}\quad}

\def\wk{well-known}
\def\socal{so-called}
\def\lhs{left-hand side}
\def\rhs{right-hand side}
\def\sfw{straightforward}

\def\CH{Cahn--Hilliard}



\def\multibold #1{\def\arg{#1}%
  \ifx\arg\pto \let\next\relax
  \else
  \def\next{\expandafter
    \def\csname #1#1#1\endcsname{{\bf #1}}%
    \multibold}%
  \fi \next}

\def\pto{.}

\def\multical #1{\def\arg{#1}%
  \ifx\arg\pto \let\next\relax
  \else
  \def\next{\expandafter
    \def\csname cal#1\endcsname{{\cal #1}}%
    \multical}%
  \fi \next}


\def\multimathop #1 {\def\arg{#1}%
  \ifx\arg\pto \let\next\relax
  \else
  \def\next{\expandafter
    \def\csname #1\endcsname{\mathop{\rm #1}\nolimits}%
    \multimathop}%
  \fi \next}

\multibold
qwertyuiopasdfghjklzxcvbnmQWERTYUIOPASDFGHJKLZXCVBNM.

\multical
QWERTYUIOPASDFGHJKLZXCVBNM.

\multimathop
ad dist div dom meas sign supp .


\def\accorpa #1#2{\eqref{#1}--\eqref{#2}}
\def\Accorpa #1#2 #3 {\gdef #1{\eqref{#2}--\eqref{#3}}%
  \wlog{}\wlog{\string #1 -> #2 - #3}\wlog{}}


\def\graffe #1{\mathopen\{#1\mathclose\}}

\def\<#1>{\mathopen\langle #1\mathclose\rangle}
\def\norma #1{\mathopen \| #1\mathclose \|}

\def\ioT {\int_0^T}
\def\intQt{\int_{Q_t}}
\def\intQ{\int_Q}
\def\iO{\int_\Omega}

\def\intQtT{\int_{Q_t^T}}

\def\inttt{\int_t^T}

\def\dt{\partial_t}
\def\dn{\partial_n}

\def\checkmmode #1{\relax\ifmmode\hbox{#1}\else{#1}\fi}
\def\aeO{\checkmmode{a.e.\ in~$\Omega$}}
\def\aeQ{\checkmmode{a.e.\ in~$Q$}}

\def\aat{\checkmmode{for a.a.~$t\in(0,T)$}}


\def\erre{{\mathbb{R}}}

\def\enne{{\mathbb{N}}}




\def\genspazio #1#2#3#4#5{#1^{#2}(#5,#4;#3)}
\def\spazio #1#2#3{\genspazio {#1}{#2}{#3}T0}

\def\L {\spazio L}
\def\H {\spazio H}

\def\C #1#2{C^{#1}([0,T];#2)}


\def\Lx #1{L^{#1}(\Omega)}
\def\Hx #1{H^{#1}(\Omega)}

\def\Ldue{\Lx 2}

\def\Huno{\Hx 1}
\def\Hdue{\Hx 2}



\let\theta\vartheta

\let\phi\varphi
\let\lam\lambda

\let\TeXchi\chi                         
\newbox\chibox
\setbox0 \hbox{\mathsurround0pt $\TeXchi$}
\setbox\chibox \hbox{\raise\dp0 \box 0 }
\def\chi{\copy\chibox}



\def\bQ{b_1}
\def\bO{b_2}
\def\bz{b_0}

\def\phQ{\phi_Q}

\def\phO{\phi_\Omega}
\def\phG{\phi_\Gamma}


\def\mG{\m_\Gamma}
\def\sO{\s_\Omega}
\def\mG{\s_\Gamma}
\def\sQ{\s_Q}

\def\Uad{\calU_{\ad}}

\def\Vp{V^*}

\def\normaV #1{\norma{#1}_V}

\let\hat\widehat


\def\cd{c_\delta}
\def\s{\sigma}  
\def\m{\mu}	    
\def\ph{\phi}	
\def\a{\alpha}	
\def\b{\beta}	
\def\d{\delta}  
\def\et{\eta}   
\def\th{\theta} 
\def\r{\rho}    
\def\g{\gamma}  
\def\bph{\bar\ph}  
\def\bm{\bar\m}    
\def\bs{\bar\s}    
\def\bu{\bar{u}}
\def\sG{\sigma^{\gamma}}    
\def\mG{\mu^{\gamma}}	    
\def\phG{\phi^{\gamma}}	    
\def\J{{\cal J}}  
\def\S{{\cal S}}   
\def\I2 #1{\int_{Q_t}|{#1}|^2}
\def\IN2 #1{\int_{Q_t}|\nabla{#1}|^2}
\def\IO2 #1{\iO |{#1(t)}|^2}
\def\INO2 #1{\iO |\nabla{#1}(t)|^2}
\def\UR{{\cal U}_R}

\def\Y{{\cal Y}}
\def\PGtil{(\tilde{P}_{\g})}

\Begin{document}


\title{Optimality conditions for an extended tumor growth
		\\[0.3cm] 
		model with double obstacle potential 
		\\[0.3cm] 
		via deep quench approach
}
\author{}
\date{}
\maketitle

\Bcenter
\vskip-1cm
{\large\sc Andrea Signori$^{(1)}$}\\
{\normalsize e-mail: {\tt andrea.signori02@universitadipavia.it}}\\[.25cm]
$^{(1)}$
{\small Dipartimento di Matematica e Applicazioni, Universit\`a di Milano--Bicocca}\\
{\small via Cozzi 55, 20125 Milano, Italy}

\Ecenter
\Begin{abstract}\noindent
In this work, we investigate a distributed optimal control problem
for an extended phase field system of Cahn--Hilliard type which physical context
is that of tumor growth dynamics.
In a previous contribution, the author has already studied the corresponding
problem for the logarithmic potential. 
Here, we try to extend the analysis by taking into account a non-smooth
singular nonlinearity, namely the double obstacle potential. 
Due to its non-smoothness behavior, the standard procedure 
to characterize the necessary conditions for the optimality cannot be performed. 
Therefore, we follow a different strategy which in
the literature is known as the ``deep quench" approach in order to obtain
some optimality conditions that
have to be interpreted in a more general framework.
We establish the existence of optimal controls and some
first-order optimality conditions for the system are derived by employing suitable
approximation schemes.

\vskip3mm
\noindent {\bf Key words:}
Distributed optimal control, tumor growth, phase field model, Cahn--Hilliard equation, 
optimal control, necessary optimality conditions, adjoint system, deep quench, asymptotic analyses.
\vskip3mm
\noindent {\bf AMS (MOS) Subject Classification:} 
35K61,  
35Q92,  
49J20,  
49K20,  
35K86,  
92C50.  

\End{abstract}

\vskip3mm

\pagestyle{myheadings}
\newcommand\testopari{\sc Signori}
\newcommand\testodispari{\sc Deep quench limit}
\markboth{\testodispari}{\testopari}

\salta
\finqui
\newpage
\section{Introduction}
\label{SEC_INTRO}
\setcounter{equation}{0}
{In recent years there has been an increased focus on the investigation and 
understanding of tumor growth by the mathematical community (see, e.g., \cite{CL}). 
Unfortunately, the phenomena which occur in real cases is far too complex to
be approached by experimental techniques as a whole. Therefore, the mathematical
modeling can be a useful instrument to reduce the problem in a manageable one since
it is able to isolate some particular mechanisms which can be hopefully sufficient to forecast 
something helpful for medical treatments.
Here, we concentrate on the subclass of models known in the literature as diffuse interface
models which are derived by continuum mixture theory.
The system we are going to consider in what follows 
constitutes a variation of the model originally 
introduced in \cite{HDZO} in order to describe the 
evolution of a young tumor, before the development of quiescent cells
in the presence of a nutrient species
(see \cite{HKNZ,WZZ,CLLW,HDPZO} and also \cite{CGH,FGR,{CGRS_ASY},CGRS_VAN}).
It consists of a \CH~equation for the phase variable
(see, e.g., \cite{Mir_CH} and the huge references therein), 
coupled with a reaction-diffusion equation for an unknown
species acting as a nutrient, for instance oxygen or glucose.
Furthermore, by interpreting the tumor and the healthy cells as inertia-less fluids,
the contribution of the velocity field can be included in the investigation by assuming
a Darcy law or a Stokes-Brinkman equation. In this regards, let us refer to
\cite{DFRGM ,GARL_1,GARL_4,GAR, EGAR, FLRS,GARL_2, GARL_3}, 
where further mechanisms such as chemotaxis and active transport are also taken into account.
Furthermore, we also point out the paper \cite{FRL}, where a non-local model is proposed.
Lastly, for different physically meaningful choices of the potentials, 
we refer to \cite{Agosti} and
to the references therein, where some numerical simulations and comparison with
clinical data can be found as well.
Further investigations and mathematical models
related to biology can be also found e.g. in \cite{DFRGM} and \cite{FLRS}.
}
Let us begin introducing our problem.
First of all, let $\Omega \subset \erre^3$ denote some open and bounded
domain with smooth boundary denoted by $\Gamma$.
For convenience, given a fixed final time $T>0$, we introduce the 
standard parabolic cylinder and its boundary by setting
\Bsist
	& \non
	Q_{t}:=\Omega \times (0,t), \quad \Sigma_{t}:=\Gamma \times (0,t)
	\quad \hbox{for every $t \in (0,T]$,}
	\\
	& Q:=Q_{T}, \quad \hbox{and} \quad \Sigma:=\Sigma_{T}.
	\label{QS}
\Esist
The problem we are going to deal with consists in minimizing the cost functional
\Bsist
	&& \non
	\J (\ph, \s, u)  := 
	\frac \bQ 2 \norma{\ph - \phQ}_{L^2(Q)}^2
	+\frac \bO 2 \norma{\ph(T)-\phO}_{L^2(\Omega)}^2
	+ \frac {b_3} 2 \norma{\s - \sQ}_{L^2(Q)}^2
	\\ && \quad
	+ \frac {b_4} 2 \norma{\s(T)-\sO}_{L^2(\Omega)}^2
    + \frac \bz 2 \norma{u}_{L^2(Q)}^2
    \label{costfunct}
\Esist
subject to the control contraints 
\Beq
    \label{Uad}
	u \in \Uad := \graffe{u \in L^{\infty}(Q): u_* \leq u \leq u^*\ \aeQ},
\Eeq
and to the state system
\begin{align}
   \a \dt \m + \dt \ph - \Delta \m &= P(\phi) (\s - \m)
  \quad \hbox{in $\, Q$}
  \label{EQprima}
  \\
   \mu &= \beta \dt \ph - \Delta \ph + \xi + \pi(\ph)
  \label{EQseconda}
  \quad \hbox{in $\,Q$}
  \\
   \xi &\in \partial I_{[-1,1]}(\ph)
  \label{EQterza}
  \quad \hbox{in $\,Q$}
  \\
   \dt \s - \Delta \s &= - P(\phi) (\s - \m) + u
  \label{EQquarta}
  \quad \hbox{in $\,Q$}
  \\
  \dn \m & = \dn \ph =\dn \s = 0
  \quad \hbox{on $\,\Sigma$}
  \label{BCEQ}
  \\
   \m(0)&=\m_0,\, \ph(0)=\ph_0,\, \s(0)=\s_0
  \quad \hbox{in $\,\Omega$.}
  \label{ICEQ}
\end{align}
\Accorpa\EQ EQprima ICEQ
At this first stage, we would confine the technicalities as much as possible, trying to 
describe the general purpose of the article. Anyhow, let us point out that
$\partial_n $ stands for the outward normal derivative, 
and $\partial I$ for the subdifferential of the indicator 
function of the interval $[-1,1]$, that is, the function which vanishes there,
and takes the value $+\infty$ elsewhere. 
Indeed, in the interval $[-1,1]$, we have the following characterization
\Beq
	\non
	s \in \partial I_{[-1,1]}(r) \quad \hbox{if and only if}  \quad 
	s \begin{cases} \leq 0 & \hbox{if } r = -1
					\\ = 0 & \hbox{if } -1 < r < 1 .
					\\ \geq 0 & \hbox{if } r = +1
    \end{cases} 
\Eeq
As far as in the above lines numerous quantities occur, let us sketch the
role they cover in our treatment, postponing the complete investigation on the
hypotheses we need to require on them to the next section. 
First we suppose $\a,\b$ to be strictly positive constant. Secondly we assume
that the symbols $\bz, \bQ,\bO,b_3,b_4$ that appear in \eqref{costfunct}
represent nonnegative constants, while $\phQ, \sQ, $ $u_*, u^*$ and $\phO,\sO $
denote some prescribed functions in their respective domains $Q$ and $\Omega$. 
Furthermore, the set $\Uad$ models the admissible set in which we can choose 
the control variables.
%

Let us spend a few words as the physical background of the
involved variables are concerned. 
The variable $\ph$, called phase variable, is
devoted to accounting for the presence of the tumor in the evolution process.
It represents a rescaled density which ranges in $[-1,1]$, where the 
extreme values $- 1$ and  $1 $ model the complete tumorous case and the healthy one, while
the values
in between denote intermediate concentrations.
The second variable $\m$ stands for the chemical potential for $\ph$,
while $\s$ is a rescaled density which takes into account the
nutrient-rich extracellular water {fraction}. Furthermore, we assume that $\s$ ranges 
from $0$ to $1$ and that when $\s$ is close to zero the nutrient is poor, 
while when $\s$ is close to one, it is rich.
Moreover, the function $P$ stands for a nonlinearity which describes the proliferation
of the tumorous cells in the tissue.
As regards $u$, it represents the control variable,
which in the application consists of a supply of a nutrient or a drug in chemotherapy.
A more complete description
of the model, along with some variations, can be found in 
\cite{CGH,FGR,{CGRS_ASY},CGRS_VAN}.
Now, let us present some literature concerning some optimal controllability results 
for these systems.
Up to our knowledge, the first control problem for a system very close to our system is 
\cite{CGRS_OPT}. There, the state system slightly differs from \EQ, since, formally,
consists of the same model in the case $\a=\b=0$. Moreover, the investigation has to be
restricted to regular potentials with polynomial growth, so that the double-obstacle one is not
allowed. In that regards, we also point out the recent \cite{CRW}, where the authors
extend the analysis of \cite{CGRS_OPT} employing a time-dependent cost functional 
which, in addition, penalizes the long-time treatments and the large mass
of the tumor at the end of the medication. 
Next, we point out \cite{S}, where the author, adding the two relaxation terms
$\a\dt\m$, $\b\dt\ph$ generalizes the optimal control problem \cite{CGRS_OPT} extending
the analysis to the case of singular and regular double-well potentials.
With the current contribution, we aim at showing that \cite{S} can be generalized allowing
also non-regular and singular potentials to be considered.
Lastly, regarding the models which take into account the velocity field, let us mention the recent contributions \cite{EK,EK_ADV}, where the authors establish the existence of optimal control and provide
some necessary and sufficient conditions for optimality (see also \cite{SW}).
 
From a heuristic viewpoint, our problem consists of choosing an admissible control variable $u$
to force the dynamics of the system to evolve towards a fixed final configuration which is 
considered to be desirable for some practical reason (e.g. surgery).

Summing up, our goal consists in solving the following problem.
\Bsist
	\non
	(P_0) && \hbox{Minimize $\J(\ph,\m,u)$}
	\hbox{ subject to the control contraints \eqref{Uad} and under the}
	\\ && \non
	\hbox{requirement that the variables $(\ph, \s)$ 
	yield a solution to \EQ.}
	\label{P_0}
\Esist
Once the well-posedness of the system \EQ~has been shown for every admissible $u$, we can define the 
corresponding control-to-state mapping $\S_0$ as the map that assigns to a fixed control
the corresponding state, namely $\S_0: u \mapsto (\m,\ph,\s)$, where $(\m,\ph,\s)$ solves \EQ.
Moreover, let us fix a convention that will be used repeatedly in the paper:
with $\S_{0,2}$ we denote the map $\S_0$ restricted to the second and 
third components, that is,
$S_{0,2}(u) = (\ph,\s)$.
Furthermore, it is possible to introduce, over $\Uad$, the 
\socal~reduced cost functional as follows
\Beq
	\label{redcost}
	\J_{\textrm{red}}(u):= \J(\S_{0,2}(u),u).
\Eeq
Thus, since we will assume $\Uad$ to be closed, bounded, convex and nonempty, 
we are able to formally characterize the necessary condition
that every optimal control $\bu$ has to satisfy by means
of the following variational inequality
\Beq
	\label{abstrnec}
	\< D{\J_{\textrm{red}}}(\bar{u}) , v - \bar{u} > \geq 0 \quad \hbox{for every } v \in \Uad,
\Eeq 
where $D{\J_{\textrm{red}}}$ denotes the derivative of \eqref{redcost} in a proper functional sense.
In this direction, it will be necessary to prove some regularity result for $\S_0$ and
we would need $I_{[-1,1]}(\cdot)$ to be differentiable, 
which is not the case. For this reason, the standard procedure 
to characterize the optimality conditions, which is essentially based on \eqref{abstrnec}, fails and
we have to proceed with a different strategy (see, e.g., \cite{Trol,Lions_OPT}). 

The technique that we are going to exploit is often referred to in the literature 
as ``deep quench" approach and it lies on a suitable approximation and
some monotonicity arguments.
As far as some recent application of this strategy is concerned,
we mention the paper {\cite{CFS}} and {\cite{CGS_DQ}}.  
The former deals with an optimal control problem for the simpler Allen--Cahn system endowed
with dynamic boundary conditions, whereas the latter
focuses the attention on the \CH~system combined with dynamic boundary conditions in which 
the control variable represents the optimal velocity for the system.
We also refer the reader to \cite{CS_doubleobstacle, CGS_nonlocal, CFGS} where,
with the same technique, other phase separation problems were faced.

Here, let us only sketch the idea behind this approximation scheme.
The key idea is the following: the differential inclusion 
\eqref{EQterza} is replaced by a function which, at every stage,
resembles the logarithmic double-well potential and which, in a proper sense, approximates 
the double obstacle nonlinearity as the parameter of the approximation
goes to zero. This plan will be quite effective since the corresponding
approximating problem complies with the framework
of {\cite{S}}, and therefore all the results there established are at our disposal.
So, we formally substitute the inclusion \eqref{EQterza}
as follows
\Beq
	\label{ideaapprox}
	\xi = g(\g)h'(\ph), 
\Eeq
where $h$ is defined by
\Beq
	\label{h}
	h(s):=
	\begin{cases} 
	(1-s)\ln(1-s) +(1+s)\ln(1+s) & \hbox{if } s\in (-1,+1)
	\\ 
	2 \ln 2 & \hbox{if } s\in \graffe{-1,+1},
	\end{cases}
\Eeq
and where $g$ is a positive and regular function acting on the parameter. 

\dafare{
We postulate that the function $g$ is positive, such that
$g \in C^0(0,1]$, and that, for every $\g \in (0,1]$, it satisfies the following
limit properties
\Bsist
	& 
	\lim_{\g \searrow 0} g(\g) = 0, 
	\quad \lim_{\g \searrow 0} \, g(\g)h(s) = 0 \quad \forall s \in[-1,+1].
	\label{g_uno}
\Esist
Moreover, the explicit expression of $h$, allows us to 
compute that $h'(s)= \ln \bigl( \frac{1+s}{1-s} \bigr)$, and
that $h''(s)= \frac 2{1-s^2} > 0 $ for $s \in (-1,1)$. Hence, combining the
requirements \eqref{g_uno} with these expressions, we also realize that
\Bsist
	\non
	\lim_{\g \searrow 0} \,g(\g)h'(s) = 0 \quad \forall s \in(-1,+1),
	\quad
	\lim_{\g \searrow 0} \,\biggl( g(\g) \lim_{s \to \pm 1} h'(s) \biggr) = \pm \infty.
\Esist	
\Accorpa\Propg h g_uno }
For instance, an admissible choice for the function $g$ could be $g(\g):= \g^p$,
for some $p>0$.
In that framework, it follows that
the graph of $g(\g)h'(\cdot)$ became closer to the one of $\partial I_{[-1,1]}(\cdot)$  as
the parameter $\g $ goes to zero. Starting from this setting, we will show that also
the family of control problems generated by this substitution, are
close to $(P_0)$ as $\g \searrow 0$ in a suitable sense.

We will see that this approximation scheme turns out to be sufficient to prove
the existence of an optimal control for the system \EQ.
Anyhow, this approach will also point out some weakness and limitations.
In particular, we will realize that nothing can be said as the 
approximation of optimal control for 
$(P_0)$ by sequences of optimal controls
for the approximating problem is concerned. This fact totally prevents
the possibility to recover some necessary conditions for the initial problem by
passing to the limit in the necessary conditions of the approximating ones.

However, something can be said if we restrict the analysis to local results.
In fact, localizing the problem around a fixed optimal control for $(P_0)$,
we can prove a kind of local density result in terms of approximating optimal controls. 
Namely, we can show that for every fixed control for $(P_0)$, say $\bu$, 
there exists an approximating sequence 
which, at every stage, is constituted by an optimal control for a suitable approximating
problem. In order to prove such an approximation result, it will be convenient to
introduce another cost functional, called adapted, which reads as follows
\Bsist
	\widetilde{\J} (\ph, \s, u):= 
	\J(\ph, \s, u)
	+ \frac 12 \norma{u-\bu}^2_{L^2(Q)}.
    \label{costfunctloc}
\Esist
Let us point out that $\widetilde{\J}$ strongly depends on the fixed control $\bu$
and also that if we take $(\ph,\s)$ 
as $(\bph, \bs)=\S(\bu)$, the adapted cost
functional $\widetilde{\J}$ reduces to $\J$. In this sense
$\widetilde{\J}$ consists of a local perturbation of $\J$ around the optimum 
$\bu$.
Next, we will solve the control problem for the approximating system subject to the adapted
cost functional and, providing to 
show some uniform estimates with respect to $\g$, we will pass to the limit in the corresponding 
optimality conditions to characterize the one for our initial system.

The plan of the rest of the paper is as follows.
In the following section, we provide a precise description of the
arguments introduce up to now at a formal level. Moreover, we fix our
setting and assumptions, and also state the established results. From the third section on, we 
focus the attention on the corresponding proofs. Furthermore, Section \ref{SEC_APPROX} is totally devoted
to the investigation of the approximating problem, while
the existence of optimal controls has been proved in Section \ref{SEC_EXOPT}.
To conclude, in Section \ref{SEC_OPTIMALITY} we perform the asymptotic analysis
that will allow us to obtain the necessary condition for the problem we are dealing with.

\section{Setting and main results}
\label{SEC_SETTING}
\setcounter{equation}{0}
In this section, we state the main results on existence of optimal
controls and on the optimality conditions.
First of all, let us define some functional spaces that will be extremely useful
later on
\Bsist
	\non
	&& H:= \Ldue, \quad V:= \Huno, \quad W:=\graffe{v \in \Hdue : \dn v = 0 \hbox{ on } \Gamma},
	\\
	\non
	&& \Y := \biggl( \H1 H \cap \L\infty V \cap \L2 W \biggr)^{\!\!3},
	\quad 
	\hat{\Y} := \Pi_2 \circ \Y ,
	\label{statespacev2}
\Esist
where $\Pi_2$ stands for the projection of the first two components. 
Furthermore, we endow them with their natural norm to obtain some Banach spaces,
and agree that the symbol $\norma{\cdot}_{X}$ denotes the norm
associated with a generic Banach space $X$. Moreover, we denote with 
$\< \cdot,\cdot>$ the duality pairing between $V$ and its dual $\Vp$.

Our basic assumptions on the system \EQ, on $\Uad$, and on the cost functional 
\eqref{costfunct} are as follows
\begin{align}
	& \bz, \bQ, \bO,b_3,b_4   \, \hbox{ are nonnegative constants, but not all zero.}
	\label{H1}
	\\ & 
	\phQ, \sQ \in L^2(Q), \phO, \sO \in \Huno, u_*,u^*\in L^\infty(Q) \hbox{ with } u_*\leq u^* \,\aeQ.
	\label{H2}
	\\ &
	\a, \b > 0.
	\label{H3}
	\\ &
 	P \in C^2(\erre) \, \hbox{ is nonnegative, bounded and Lipschitz continuous.}
	\label{H4}
	\\ &
	\hat{\pi} \in C^3(\erre) \,\, \hbox{and} \,\, \pi:=\hat{\pi}' \hbox{ is	\Lip~continuous.}
	\label{H5}
	\\ &
	\m_0 , \ph_0 , \s_0 \in \Huno.
	\label{H6}
	\\ &
	{|\ph_0| \leq 1 \quad a.e. \hbox{ in } \Omega.}
	\label{H7}
\end{align}
We thus may infer that the choices of $u_*$ and $u^*$ entails that the set $\Uad$
in turn is bounded, closed, convex and nonempty.
Moreover, given a positive constant $R$, we introduce
\begin{align*}
	\non &
	\UR \subset L^2(Q) \hbox{ be a nonempty and bounded open set such that
	it contains } \Uad, 
	\\ \non & 
	\hbox{and } \norma u_2 \leq R \hbox{ for all }
	u \in \UR.
\end{align*}
{Let us emphasize that the singular, while not regular,
double-well potential we are considering is usually referred to as the double-obstacle 
potential and it reads as
\Beq
	\non
	F_{2obst} : = I_{[-1,1]} + \hat{\pi}.
\Eeq
One has also to keep in mind that  $\partial I_{[-1,1]} $  may be multivalued and therefore,
we introduce a selection $\xi$ by \eqref{EQterza} that
may not be regular enough to possess a trace.}

Now, we can start to list our results.
\Bthm[Well-posedness]
\label{THM_Existenceeq}
Suppose that \eqref{H3}-\eqref{H7} are fullfilled. Then, for every $u\in L^2(Q)$, there exists a
unique quadruplet $(\m,\ph,\s,\xi) \in \Y \times L^2(Q)$ which solves \EQ.
\Ethm
\proof[Proof of Theorem \ref{THM_Existenceeq}]
It directly follows as a special case from {\cite[Thm.~2.2, p.~2426]{CGH}} and 
{\cite[Thm.~2.2, p.~97]{CGRS_VAN}}.
$\qed$

The above result allows us to properly introduce
the control-to-state mapping $\S_0: L^2(Q)\to \Y$, which is
the map that assigns to every admissible control $u$ the corresponding solution
$(\m,\ph,\s)$ to system \EQ. 
As the approximating system is concerned,
assuming \Propg, and replacing \eqref{EQterza} by \eqref{ideaapprox}, we
get the following system which depends on $\g$ and reads as
\begin{align}
   \a \dt \mG + \dt \phG - \Delta \mG &= P(\phG) (\sG - \mG)
  \quad \hbox{in $\, Q$}
  \label{EQGprima}
  \\
   \mG &= \beta \dt \phG - \Delta \phG + g(\g)h'(\phG) + \pi(\phG)
  \label{EQGseconda}
  \quad \hbox{in $\,Q$}
  \\
   \dt \sG - \Delta \sG &= - P(\phG) (\sG - \mG) + u
  \label{EQGterza}
  \quad \hbox{in $\,Q$}
  \\
  \dn \mG & = \dn \phG =\dn \sG = 0
  \quad \hbox{on $\,\Sigma$}
  \label{BCEQG}
  \\
   \mG(0)&=\m_0^{\gamma},\, \phG(0)=\ph_0^{\gamma},\, \sG(0)=\s_0^{\gamma}
  \quad \hbox{in $\,\Omega$},
  \label{ICEQG}
\end{align}
\Accorpa\EQG EQGprima ICEQG
where $\graffe {(\mG_0,\phG_0,\sG_0)}_\g$ denotes a family of initial data.
We postulate that such a family fulfills the following requirements
\Bsist
	\label{dataapproxuno}
	& (\mG_0,\phG_0,\sG_0) \in 
			(\Huno \cap L^\infty(\Omega)) \times W \times \Huno 
			\quad \forall \g\in(0,1],
	\\  \label{dataapproxdue}
	& |\phG_0| \leq 1- \g/2 \quad \aeO   \quad \forall \g\in(0,1],
	\\ \label{dataapproxtre}
	& (\mG_0,\phG_0,\sG_0)  \to (\m_0,\ph_0,\s_0) \quad \hbox{strongly in } V \times V \times V 
			\quad \hbox{as } \g\searrow 0.
\Esist
\Accorpa\DataApprox dataapproxuno dataapproxtre
Even though they seem to be reasonable assumptions, the existence of a family that satisfies all the 
above conditions is not so trivial. For this reason we refer to Appendix \ref{APPENDIX}, 
where the construction of such a family is shown in detail.
Since the approximating problem \EQG~perfectly fits the framework of
\cite{S}, the result below follows by a simple application of {\cite[Thm.~2.1]{S}}.

\Bthm[Well-posedness of the approximating system]
\label{THM_Gexistence}
Assume that \Propg, \eqref{H3}-\eqref{H6}, and \DataApprox~are in force. 
\begin{enumerate}
	\item[(i)] For every $\g\in (0,1]$ and every $u\in\Uad$, the system \EQG~admits a unique solution 
	$(\mG,\phG,\sG)$ that possesses the following regularity
	\Bsist
	\label{reg_gamma_ph}
	& \phG \in W^{1,\infty}(0,T;H) \cap \H1 V \cap \L\infty W \subset \C0 {C^0({\overline{\Omega}})}
	\\ 
	\label{reg_gamma_ms}
	&\mG, \sG \in \H1 H  \cap \L\infty V \cap \L2 W \subset C^0([0,T];V)
	\\ 
	\label{reg_gamma_minf}
	&\mG \in L^{\infty}(Q),
	\Esist
	and whose second component also satisfies that $$ -1<\inf\phG \leq \sup\phG <1 \quad \aeQ. $$
	\item[(ii)]For every given $\g\in(0,1]$ there exist constants 
	$\ph_*(\g),\ph^*(\g)\in(-1,1)$, which depend on $\g$, on the initial data, 
	and on the data of the system such that
	\Beq
		\label{sep_g}
		\ph_*(\g)\leq \phG \leq \ph^*(\g)\quad \aeQ,
	\Eeq
	where $\phG$ is the second component of the unique solution to the approximating 
	system \EQG~associated to the given $u$.
\end{enumerate}
\Ethm

\Brem
Let us point out that, unfortunately, a uniform separation property from $(-1,1)$ is out of reach.
In fact, when $\g\searrow 0$, it may occur that
\Beq
	\non
	\ph_*(\g)\searrow -1 \quad \hbox{and/or} \quad \ph^*(\g)\nearrow 1.
\Eeq
\Erem

Nonetheless, although we cannot prove a separation result for
\EQ, we will see that in the asymptotic investigation the separation 
property \eqref{sep_g} for the approximated system will be using several times.
Another ingredient that will be fundamental in the asymptotic analysis is
the result below.

\Blem
\label{LM_regG}
Suppose that $u\in\UR$, and that assumptions \Propg, \eqref{H3}-\eqref{H6}, and \DataApprox~are fulfilled. 
Then, whenever $\g\in(0,1]$ and $(\mG,\phG,\sG)$ is the corresponding 
solution to \EQG, it holds that
\Beq
	\label{norma_Y} 
	\norma{(\mG,\phG,\sG)}_\Y \leq C_2,
\Eeq
where $C_2$ is a positive constant which only depends on the data of the system, 
but it is independent of the parameter $\g$. Moreover, 
we have the following estimate
\Beq
	\label{norma_gh}
	\norma{g(\g)h'(\phG)}_{L^2(Q)}\leq C_3.
\Eeq
\Elem

Once the state system \EQG~has been analyzed, we can address the corresponding control
problem and, as above, we set the following minimization problem.
\Bsist
	\non
	(P_\g) && \hbox{Minimize $\J(\phG, \sG, u)$}
	\hbox{ subject to the control contraints \eqref{Uad} and under }
	\\ && \non
	\hbox{the requirement that the variables $(\phG, \sG)$ are the components of the }
	\\ && \non
	\hbox{solution to \EQG.}
	\label{P_G}
\Esist

{Let us recall that 
in Theorem \ref{THM_Gexistence} the well-posedness of system \EQG~has been already shown.
Therefore, we are also in a position to define, for every $\g\in(0,1]$, 
the well-posed map $\S_\g$ consisting of the control-to-state mapping associated
to the system \EQG, and the corresponding restriction $\S_{\g,2}$.
In addition, since $(P_\g)$ complies with the framework of {\cite{S}},
a simple application of {\cite[Thm.~2.6]{S}}
leads to the following lemma.
\Blem
\label{LM_existencePG}
Assume that \Propg, \eqref{H1}-\eqref{H7}, and \DataApprox are in force. 
Then, whenever $\g\in(0,1]$ is given, the problem $(P_\g)$ admits at least
an optimal control.
\Elem
}
Now, let us present the existence result, whose proof will be structured
this way: first, we approximate in a suitable sense $(P_0)$ as $\g \searrow 0$ by $(P_\g)$ 
and then, accounting for some compactness and monotonicity arguments, 
we pass to the limit.
\Bthm[Existence of optimal controls]
\label{THM_Existenceopt}
Suppose that \eqref{H1}-\eqref{H7} are fulfilled. Then $(P_0)$ admits at least a solution.
\Ethm
As already mentioned, the existence result combined with some asymptotic techniques 
turns out to be insufficient to properly characterize the optimality conditions
we are looking for. In particular, the fact that every optimal control 
for $(P_0)$ can be found as a limit, in a 
proper topology, of some approximating sequences of optimal controls for $(P_{\g_n})$
cannot be proven, whenever $\graffe{\g_n}\subset (0,1]$ denotes a sequence which goes to zero as 
$n \nearrow \infty$. 
This gives no hope to recover some necessary conditions for $(P_0)$ by
the mere investigation of the ones of $(P_{\g_n})$, since the passage
to the limit at that stage will be meaningless.

{However, even though we are not able to prove such a global result, a partial one can be
stated localizing the problem.
The key idea, which was introduced by Barbu in {\cite{BARBU}}, consists of
investigating the same approximating problem \EQG, but focusing the attention
on the control problem corresponding to the adapted cost functional.}
Hence, the \socal~adapted control problem reads as follows:
\Bsist
	\non
	(\widetilde{P}_\g) && \hbox{Minimize $\widetilde{\J}(\phG, \sG,u)$}
	\hbox{ subject to the control contraints \eqref{Uad} and under the}
	\\ && \non
	\hbox{requirement that the variables $(\phG, \sG)$ yield a solution to \EQG,}
	\label{P_Gtilda}
\Esist
where let us remind that $\widetilde{\J}$ is defined, once that $\bu$ has been fixed, by 
$$\widetilde{\J} (\ph, \s, u):= \J(\ph, \s, u) + \frac 12 \norma{u-\bu}^2_{L^2(Q)}.$$
Again, we can easily prove the following lemma which \sfw ly~follows
as an application of {\cite[Thm.~2.6]{S}}.
\Blem
\label{LM_existencePGtil}
Under the assumptions \Propg, \eqref{H1}-\eqref{H6}, and \DataApprox,
whenever $\bar{u}\in \Uad$ and $\g\in(0,1]$ are given,
the optimal control problem $(\widetilde{P}_\g)$ possesses at least a solution.
\Elem
The key result which motivate the interest toward the adapted optimal control 
problem is formalized in 
the next theorem where we will show that every fixed optimal control for $(P_0)$
can be obtained as a limit of a sequence of optimal controls
for $(\widetilde{P}_\g)$ which is of great importance for the forthcoming asymptotic analysis. 
Indeed, we have:
\Bthm
\label{THM_Approximation}
Assume that \Propg, \eqref{H1}-\eqref{H7}, and \DataApprox~are in force. 
Moreover, let $(\bph,\bs,\bu)\in\hat{\Y}\times\Uad$ be an
optimal choice for $(P_0)$. Then, for every sequence $\graffe{\g_n}_n$ which goes
to zero as $n\nearrow\infty$, and for every $n \in \enne$, there exists an 
approximating optimal triple for $(\widetilde{P}_{\g_n})$, namely 
$(\bph^{\g_n},\bs^{\g_n},\bu^{\g_n})$, whose 
components satisfy the following convergences as $n\nearrow\infty$
\Bsist
	&& \label{stimaadapteduno}
	\bu^{\g_n} \to \bu \ \ \hbox{strongly in } L^2(Q)
	\\  && \label{stimaadapteddue}
	\bph^{\g_n} \to \bph \ \ \hbox{weakly star in }  \H1 H \cap \L\infty V \cap \L2 W
	\\  && \label{stimaadaptedtre}
	\bs^{\g_n} \to \bs \ \ \hbox{weakly star in }  \H1 H \cap \L\infty V \cap \L2 W
	\\  && \label{stimaadaptedquattro}
	\widetilde{\J} (\bph^{\g_n},\bs^{\g_n},\bu^{\g_n}) \to \J(\bph,\bs,\bu).
\Esist
\Accorpa\Conv stimaadapteduno stimaadaptedquattro
\Ethm
Let us emphasize that this is the correct way in which the assertion 
``$(P_\g)$ approximates $(P_0)$ as $\g\searrow0$" has to be interpreted.
By virtue of the above result, it is somehow reasonable that some 
optimality conditions can be earn
exploiting this result. Indeed, we will establish 
the first optimality conditions for $(P_0)$ by passing to the limit
as $n \nearrow \infty$ in the corresponding optimality conditions 
for the approximating problem $(\widetilde{P}_{\g_n})$.

In Section \ref{SEC_OPTIMALITY}, we will present the adjoint system for 
\EQG~which was originally treated in \cite{S}. That system is formulated in terms of the variables 
${q}^{\g},{p}^{\g},{r}^{\g}$ and it represents a core argument
in order to pass to the limit to characterize the optimality conditions for \EQ.
At this stage, let us only disclose that it admits existence and uniqueness 
of a solution and also that we are able to show the regularity that its solution
enjoys (cf. Theorem \ref{THMadjoint}).
Next, we investigate the properties of this system
in order to let $\g\searrow 0$ to characterize the optimality conditions
we are looking for.
To precisely state the asymptotic 
result we have established, let us first introduce further spaces that
will naturally come out from the mathematical analysis. We set
\Bsist
	\non
	& {\cal Z} := (\L\infty H \cap \L2 V) \times (\H1 H \cap \L\infty V \cap \L2 W)
	\\ & \quad  
	\times (\H1 H \cap \L\infty V \cap \L2 W),
	\\ 
	&{\cal W}(0,T) := \H1 \Vp \cap \L2 V \subset \C0 H,
	\\
	&{\cal W}_0 (0,T) := \graffe{ v \in {\cal W}(0,T) : v(0)=0},
\Esist
and endow these spaces with their natural norm to get three Banach spaces.
To avoid an heavy notation, we will denote the norm of ${\cal W} (0,T)$ and ${\cal W}_0 (0,T)$
by $\norma{\cdot}_{{\cal W} }$ and $\norma{\cdot}_{{\cal W}_0}$, respectively.
Furthermore, we convey to use ${}_{{\cal W}_0^* (0,T)}\< \cdot , \cdot >_{{\cal W}_0(0,T)}$
for the duality product between the dual of ${\cal W}_0 (0,T)$, 
${\cal W}_0 (0,T)^*$ and ${\cal W}_0 (0,T)$ itself.
Moreover, it is worth noting that the space $\L2 \Vp$ is embedded into ${\cal W}_0 (0,T)^*$.
In fact, if $z$ belongs to $\L2 \Vp$, for every $v$ in ${\cal W}_0 (0,T)$, we have that 
\Beq
	\non
	\quad_{{\cal W}_0 (0,T)^*} \< z, v>_{{\cal W}_0 (0,T)}
	= \ioT {}_{\Vp}\< z(t) , v (t) >_{V} \,dt.
\Eeq

In conclusion, we present the results we are able to prove which specify
how the new optimality conditions have to be read.
\Bthm
\label{THM_stimeadj}
Suppose that \Propg, \eqref{H1}-\eqref{H7}, and \DataApprox~are satisfied, 
and let us define
\Beq
	\label{deflambda}
	\lam^{\g_n} := g(\g_n)h''(\bph^{\g_n})\, \bar{q}^{\g_n} \quad \hbox{ for every } n \in \enne,
\Eeq
where $h$ is the function introduced by \eqref{h} and where $\graffe{\g_n} \subset (0,1]$ 
is a sequence which goes to zero as $n \nearrow \infty$.
Then, there exists a positive constant $C_4$ such that, for every $n \in \enne$,
the following holds true
\Beq
	\label{stimaadj}
	\norma{(\bar{q}^{\g_n},\bar{p}^{\g_n},\bar{r}^{\g_n})}_{\cal Z} 
	+ \norma{\dt \bar{q}^{\g_n}}_{{\cal W} _0(0,T)^*} 
	+ \norma{\lam^{\g_n}}_{{\cal W}_0 (0,T)^*} 
	\leq C_4,
\Eeq
where the variables $\bar{q}^{\g_n},\bar{p}^{\g_n},\bar{r}^{\g_n}$ denote the unique 
solutions to the adjoint problem for \EQG~considered for the element of the 
sequence $\g_n$ instead of for $\g$,
and where the constant $C_4$ may depend on the data of the system, but it is independent of $n$.
In addition, up to a subsequence, we deduce the following convergences
\Bsist
	&& \label{conv_q}
	\bar{q}^{\g_n} \to q\ \ \hbox{weakly star in } \L\infty H \cap \L2 V
	\\ &&
	\label{conv_p}
	\bar{p}^{\g_n} \to p \ \ \hbox{weakly star in } \H1 H \cap \L\infty V \cap \L2 W
	\\  && \label{conv_r}
	\bar{r}^{\g_n} \to r\ \ \hbox{weakly star in } \H1 H \cap \L\infty V \cap \L2 W
	 \\  && \label{conv_lam}
	 \lam^\g \to \lam \ \ \hbox{weakly in } {{\cal W}_0^*}.
\Esist
\Accorpa\Conv conv_q conv_lam
\Ethm
This will allow us to let $\g\searrow 0$ and prove the following optimality conditions.
\Bthm
\label{THM_optimality}
Let the assumptions \accorpa{h}{g_uno}, \eqref{H1}-\eqref{H7}, and \DataApprox~be 
fulfilled. Furthermore, let
$(\bph,\bs,\bu)\in\hat{\Y}\times\Uad$ be an optimal choice for $(P_0)$.
Then the following properties hold true.\\
(i) Whenever a sequence $\graffe {\g_n}\subset (0,1]$ which goes to zero as 
$n\nearrow\infty$ is fixed, we have that for every $n\in\enne$ there exists an approximating
optimal triple $(\bph^{\g_n},\bs^{\g_n},\bar{u}^{\g_n})$ 
to $(\bph,\bs,\bu)$, namely a triplet which solves the adapted control problem 
$(\widetilde{P}_{\g_n})$ and which, as $n\nearrow\infty$, satisfies the convergences pointed out by \Conv.\\
(ii) Moreover, under the same assumptions we have that for every subsequence 
$\graffe {n_k}_k$ of $\enne$, there exists a subsequence $\graffe {n_k{_j}}_j$, 
a triple $(p,q,r)\in \cal Z$, and a functional $\lam \in {\cal W}_0(0,T)^*$ such that 
the variational inequality which characterizes the optimality
\Beq
	\label{nec_final}
	\intQ ({r} + \bz \bu)(v - \bu)
	\geq 0 \quad \forall v \in \Uad,
\Eeq
is satisfied. Furthermore, \aat, the triplet $(q,p,r)$ solves the adjoint system 
consisting of the following variational equation
\begin{align}
	\non
	&  \ioT \<\dt v , p -\b q> 
	- {}_{{\cal W}_0^*}\< \lam, v >_{{\cal W}_0}
	+ \intQ \nabla q \cdot \nabla v
	- \intQ \pi''(\bph) q v
	\\ & \quad
	+ \intQ P'(\bph)(\bs - \bm)(r-p)\, v
	= 
	\intQ \bQ (\bph - \phQ) v
	+ \iO \bO(\bph(T) - \phO) v(T)
	,
	\label{adj_FV_uno}
\end{align}
which holds for every $v\in {\cal W}_0 (0,T)$, combined with
\begin{align}
	\quad
	q -\a \dt p - \Delta p 
	+ P(\bph)(p -r) &=   0
    \quad \hbox{in $\,Q$}
    \label{adj_FV_due}
    \\  
    -\dt r - \Delta r + P(\bph)(r - p) 
    &= {b_3}(\bs - \sQ)
    \quad \hbox{in $\,Q$}
    \label{adj_FV_tre}
    \\ 
    \dn p &= \dn r = 0
    \quad \hbox{on $\Sigma$}
    \label{adj_FV_boundary}
    \\  
	\a p(T) &= 0, 
	\,\, r(T) = b_4 (\bs(T) - \sO)
	\quad \hbox{in $\Omega$.}
	\label{adj_FV_fin}
\end{align}
\Accorpa\ADJfin adj_FV_uno adj_FV_fin
In addition, we can show that 
\Bsist
	\label{slackness_uno}
	\liminf_{n\to\infty} \intQ \lam^{\g_n} \bar{q}^{\g_n} \geq 0,
\Esist
and also a complementary slackness condition 
\Bsist
	\label{slackness_due}
	\lim_{n\to\infty} \intQ \lam^{\g_n}(1-(\bph^{\g_n})^2) \Phi
	= 2\, \lim_{n\to\infty} \intQ g(\g_n) \, \bar{q}^{\g_n}\Phi =0,
\Esist
are satisfied, where $\Phi$ is a general smooth function which vanishes at zero.
Note that the above limits should be considered, in principle, for the 
subsequence $\graffe {n_k{_j}}_j$ as $j\nearrow\infty$.
\Ethm
\Brem 
It can possibly occur that the limit triple $(q,p,r)$, and/or $\lam$ may be not uniquely
determined. Indeed, it should happen that for different subsequences the corresponding limit 
$q,p,r$, and $\lam$ may change. Anyhow, a proper kind of uniqueness can be stated
in term of a suitable projection. 
In fact, whenever $\bz>0$, it follows from the variational inequality
\eqref{nec_final} that
\Beq
	\non
	\bu = \mathbb{P}_{\Uad}(-\bz^{-1}r) ,
\Eeq
where $\mathbb{P}_{\Uad}$ stands for the orthogonal projection onto
${\Uad}$ with respect to the standard inner product of $L^2(Q)$, so
that $r$ is uniquely determined in terms of $\bu$.
\Erem
In the remainder, we recollect some useful inequalities and properties that are
applied several times in the paper.
First, we often owe to the \wk~Young inequality
\Beq
  ab \leq \delta a^2 + \frac 1 {4\delta} \, b^2
  \quad \hbox{for every $a,b\geq 0$ and $\delta>0$}.
  \label{young}
\Eeq
Furthermore, since the evolution set $\Omega$ is a bounded subset of
$\erre^3$ and possesses regular boundary, we can account for the 
Sobolev continuous embedding
\Beq
	\label{sobolev_emb}
	\Huno \hookrightarrow L^q(\Omega)
	\quad \hbox{which is satisfied for every $q\in [1,6]$,}
\Eeq
i.e. we have the existence of suitable constant for which the following
inequality holds true
\Beq
	\non
	\norma{v}_{L^q(\Omega)} \leq C_q \norma{v}_{V}	
	\quad \hbox{for every $v \in V$ and $q\in [1,6]$.}
\Eeq
\Brem
\label{convcostanti}
Let us conclude the section explaining a convention that we are going to employ
as far as constants are concerned. Since in the following we have to 
deal with several estimates, we agree to use the symbol $c$ for every constant
which depends only on the final time~$T$, on~$\Omega$, the shape of the nonlinearities,
on the norms of the involved functions, and possibly on 
$\a$ and $\b$, but it has to be independent of $\g$. Therefore, the meaning of $c$ might
change from line to line and even in the same chain of inequalities.
On the other hand, the capital letters are devoted to denoting precise constants
which we eventually will refer to.
\Erem

\section{Approximating system}
\label{SEC_APPROX}
\setcounter{equation}{0}
From this section on, we start with the proofs of the stated results.
In this section, we deal with the approximating system \EQG.
Since a lot of properties immediately follows from {\cite{S}}
no repetition is required here. As a matter of fact, we only need to prove Lemma \ref{LM_regG},
namely check that the constant $C_2$ involved in the lemma
turns out to be independent of $\g$. 
This property will be fundamental later on to let $\g \searrow 0$ in order to 
to prove the existence of an optimal control
for $(P_0)$.

\proof[Proof of Lemma \ref{LM_regG}]
We employ similar estimates to the ones performed 
in {\cite[Proof of Thm.~2.1]{S}} while referred to the 
approximating problem \EQG. Furthermore, we will have
the care to show that all the constants that will appear do not depend on $\g$.

{\bf First estimate:} 
We add to both sides of \eqref{EQGseconda} the term $\phG$,
multiply \eqref{EQGprima} by $\mG$, this new second equation by $-\dt\phG$ and \eqref{EQGterza}
by $\sG	$, then we add the resulting equations and integrate over $Q_t$ and by parts
to obtain
\Bsist
	\non &&
	\frac \a2 \IO2 \mG
	+ \I2 {\nabla\mG}
	+ \b \I2 {\dt\phG}
	+ \frac 12 \IO2 {\phG}
	+ \frac 12 \IO2 {\nabla\phG}
	\\ \non && \quad
	+ \,g(\g)\iO h(\phG(t))
	+ \frac 12 \IO2 \sG
	+ \I2 {\nabla\sG}
	+ \intQt P(\phG)(\sG-\mG)^2
	\\ \non && \quad
	\leq 
	\frac \a2 \iO |\mG_0|^2
	+ \frac 12 \iO |\phG_0|^2
	+ \frac 12 \iO |\nabla\phG_0|^2
	+  g(\g)\iO h(\phG_0)
	\\ \non && \quad
	+ \frac 12 \iO |\sG_0|^2
	+ \intQt u\sG
	+ \intQt \phG\,\dt\phG
	+ \intQt \pi(\phG) \dt \phG,
\Esist
where almost all the integrals on the \lhs~are nonnegative since they are 
squares and $P$ is so by \eqref{H4}.
In fact, the only term that needs further manipulations is the sixth,
the one in which $g$ and $h$ appear. 
{First of all, let us recall that $g$ is assumed to be positive,
that $h$ is defined by \eqref{h} and that it 
remains bounded in the interval $(-1,1)$. Moreover, as a solution, $\phG$ 
possesses the regularity stated by \eqref{reg_gamma_ph}
and also enjoys the separation result \eqref{sep_g}.
Therefore, both $\phG(t)$ and $\phG_0$ range in the inner part of the interval
$(-1,1)$, and from this property, along with the fact that $h$ is bounded from below in $[-1,1]$,
we infer that $h(\phG(t))$ is bounded from below as well.
Hence, we neglect that term and we are reduced to control 
the integrals on the \rhs~which we denote
by $I_1,...,I_8$, in this order. 
Due to the requirements \DataApprox~and the above observation we immediately
deduce that 
\Beq
	\non
	\sum_{i=1}^5 |I_i| 
	\leq
	c.
\Eeq 
As the other terms are concerned, we invoke the Young inequality to show that }
\Bsist
	\non
	|I_6|+|I_7|+|I_8| 
	\leq
	\frac 12 \I2 u
	+ \frac 12 \I2 {\s^\g}
	+ 2\d \I2 {\dt {\ph^\g}}
	+ \cd \intQt ({|\ph^\g|^2}+1),
\Esist
for a positive $\d$ yet to be determined.
Hence, we fix $0<\d<\b/2$ and a Gronwall argument yields
\Bsist
 	\non &&
 	\norma{\mG}_{\L\infty H \cap \L2 V}
	+ \norma{\phG}_{\H1 H \cap \L\infty V}
	+ \norma{\sG}_{\L\infty H \cap \L2 V}
	\leq 
	c.
	\label{Iest}
\Esist

{\bf Second estimate:} 
Now, we multiply \eqref{EQGprima} by $\dt\mG $ and \eqref{EQGterza}
by $\dt\sG$, add the resulting equalities and integrate over $Q_t$.
Owing to the above estimate we easily conclude that
\Bsist
	\non
	&& \norma{\mG}_{\H1 H \cap \L\infty V}
	+ \norma{\sG}_{\H1 H \cap \L\infty V}  
	\leq
	c.
	\label{IIest}
\Esist

{\bf Third estimate:} 
Equations \eqref{EQGprima} and \eqref{EQGterza} show a parabolic structure with respect to
the variables $\mG$ and $\sG$, respectively. In turn, it follows from the above estimates
that their forcing terms both belong to $\L2 H$.
Hence, since the initial data \eqref{ICEQG} 
are regular due to \eqref{dataapproxuno}, the
{elliptic regularity theory for homogeneous Neumann boundary problems produces}
\Bsist
	\non
	\norma{\mG}_{\L2 W}
	+ \norma{\sG}_{\L2 W}
	\leq 
	c.
	\label{IIIest}
\Esist

{\bf Fourth estimate:}
Next, we aim at improving the regularity of $\phG$ arguing in a similar way 
via a comparison argument in \eqref{EQGseconda}. 
So, let us rearrange \eqref{EQGseconda} as follows
\Bsist
	\label{fortheq}
	-\Delta \phG + g(\g)h'(\phG) = f, \quad \hbox{where $f:=\mG -\b \,\dt \phG -\pi(\phG)$,}
\Esist
and we realize that the above estimates entail that $f \in \L2 H$.
We then test \eqref{fortheq} by $- \Delta \phG$ and integrate over $\Omega$ to get, 
\aat, the following inequality
\Bsist
	\non
	\IO2 {\Delta\phG}
	+ g(\g) \iO h''(\phG(t)) \, |\nabla \phG(t)|^2 
	\leq 
	\frac 12 \IO2 {\Delta\phG}
	+\frac 12 \IO2 {f},
\Esist
where we have applied Young's inequality on the \rhs.
The second term of the \lhs~turns out to be positive since $\phG(t)$
satisfies \eqref{sep_g} and $h''$ is nonnegative in such an interval. 
Hence, we realize that
\Beq
	\non
	\frac 12 \IO2 {\Delta\phG}
	\leq
	\frac 12 \IO2 {f},
\Eeq
and the elliptic regularity theory for homogeneous Neumann problem yields that
\Beq
	\non
	\norma{\phG}_{\L2 W} 
	\leq c.
	\label{IVest}
\Eeq
Combining all the above estimates, it is \sfw~to realize that \eqref{norma_Y} has been proved.

Let us conclude by proving the second part of the lemma.

{\bf Fifth estimate:}
On account of the above estimates, a comparison in \eqref{EQGseconda} directly gives that
\Beq
	\label{stima_xi}
	\norma{g(\g)h'(\ph^\g)}_{L^2(Q)} \leq c,
\Eeq
that is the estimate we are looking for.
\qed

\section{Existence and approximation of optimal controls}
\label{SEC_EXOPT}
\setcounter{equation}{0}
Here, we essentially aim to show the validity of Theorems \ref{THM_Existenceopt} and
\ref{THM_Approximation}. 

\proof[Proof of Theorem \ref{THM_Existenceopt}]
Let us pick an arbitrary sequence $\graffe {\g_n}\subset(0,1]$ which goes to zero as
$n \nearrow \infty$. In view of Lemma \ref{LM_existencePG}, 
we can take an optimal triple for $(P_{\g_n})$ 
associated with that sequence. 
Namely, for every $n\in\enne$, we consider the following triple
\Beq
	(\ph^{\g_n},\s^{\g_n},u^{\g_n})\in \hat{\Y} \times \Uad,
\Eeq
where $(\ph^{\g_n},\s^{\g_n})= \S_{\g_n,2}(u^{\g_n})$.
Moreover, from Lemma \ref{LM_regG} and Theorem \ref{THM_Gexistence} it follows that
\Bsist
	|\ph^{\g_n}|<1 \quad \aeQ, \aand \norma{(\ph^{\g_n},\s^{\g_n})}_{\hat{\Y}} \leq C_2
	\quad \hbox{for every $n\in\enne$,}
\Esist
where the constant $C_2$ is independent of $\g$.
Therefore, thanks to \wk~weak star compactness, it is a 
standard matter to show that there exist some $\bu\in\Uad$, and a 
triple $(\bm,\bph,\bs)\in {\Y}$, such that the following convergences
\Bsist
	&& \non
	u^{\g_n} \to \bu \ \ \hbox{weakly star in } L^\infty(Q)
	\\  && \non
	\ph^{\g_n} \to \bph \ \hbox{weakly star in }  \H1 H \cap \L\infty V \cap \L2 W
	\\  && \non
	\s^{\g_n} \to \bs \ \ \hbox{weakly star in }  \H1 H \cap \L\infty V \cap \L2 W
	\\ && \non
	\m^{\g_n} \to \bm \ \ \hbox{weakly star in }  \H1 H \cap \L\infty V \cap \L2 W,
\Esist
hold as $n \nearrow \infty$ .
Furthermore, due to the continuity of the embedding
\Beq
	\non
	\H1 H \cap \L2 W \subset \C0 V,
\Eeq
we infer that $\bph,\bs \in \C0 V$.
In addition, also some strong convergences can be recovered invoking 
the Aubin-Lions lemma (see, e.g., \cite[Sect.~8, Cor.~4]{Simon}). 
Indeed, one can show that
\Bsist
	\label{strongconvuno}
	&&\ph^{\g_n} \to \bph \ \ \hbox{strongly in }  \C0 H \cap \L2 V 
	\\  && \label{strongconvdue}
	\s^{\g_n} \to \bs \ \ \hbox{strongly in }  \C0 H \cap \L2 V
	\\  && \label{strongconvtre}
	\m^{\g_n} \to \bm \ \ \hbox{strongly in }  \C0 H \cap \L2 V.
\Esist
At this point, from the assumptions on the initial data \DataApprox, 
combined with the above strong convergences, we infer
that $\bph(0)=\ph_0$ and $\bs(0)=\s_0$. By the same token, we can handle both the 
nonlinearities $P$ and $\pi$, that are \Lip~continuous by \eqref{H4} 
and \eqref{H6}, respectively. In fact, as $n \nearrow \infty$, we also realize that
\begin{align}
	\non 
	P({\ph^{\g_n}}) &\to P(\bph) \ \ \hbox{strongly in } \L2 H
	\\ \non  
	\pi{(\ph^{\g_n})} &\to \pi(\bph)\ \ \hbox{strongly in }  \L2 H.
\end{align}
Moreover, estimate \eqref{stima_xi} also leads to the following weak convergence 
\Beq
	\label{weak_conv}
	g(\g_n)\,h'(\ph^{\g_n}) \to \xi \ \ \hbox{weakly in } \L2 H.
\Eeq
Now, we claim that its limit, that we have denoted by $\xi$, has the same meaning of the 
variable $\xi$ introduced in \eqref{EQseconda},
namely that $\xi\in\partial I_{[-1,1]}(\bph) \,\, \aeQ$. 
For this purpose, we account for the  convexity property of $h$ defined by \eqref{h}. 
It yields that, for every $n\in\enne$, the following inequality 
is satisfied
\Bsist
	\non
	\intQ g(\g_n)\, h(\ph^{\g_n})
	+\intQ g(\g_n)\,h'(\ph^{\g_n})(w-\ph^{\g_n})
	\leq 
	\intQ g(\g_n)\,h(w)
	\,\,\, \hbox{for all $ w \in B_1$,}
\Esist
where $B_1:=\graffe{ v \in L^2 (Q): |v| \leq 1 \,\, \aeQ}$.
Then, we owe to \eqref{g_uno} and deduce that in the above inequality 
the former and the latter terms go to zero as $n \nearrow \infty$ .
Therefore, we combine the strong convergence \eqref{strongconvuno} 
with the weak one \eqref{weak_conv} to find in the limit
\Beq
	\intQ \xi (\bph - w)\geq 0 \quad \hbox{ for every $w\in B_1$,}
\Eeq
which means exactly that $\xi$ is an
element of the subdifferential of the
extension of $I_{[-1,1]}(\cdot)$ to $L^2(Q)$, or equivalently 
(cf. \cite[Ex.~2.3.3., p.~25]{BRZ}) that 
$\xi \in \partial I_{[-1,1]}(\bph)$, as we claimed.
Hence, all the above convergences ensure us the possibility to pass to the limit
as $n\nearrow\infty$ and realize that $(\bph,\bs,\bu)$ solves \EQ~and
also that it is an admissible choice for $(P_0)$, i.e. that $(\bph,\bs)=\S_{0,2}(\bu)$.

To conclude, it remains to show that $(\bph,\bs,\bu)$ is in fact not only admissible, but
also optimal and we will manage this problem by accounting for monotonicity arguments.
In fact, recalling that \eqref{costfunct} is lower semicontinuous, it turns out that
the following
\begin{align}
	& \non
	\J(\bph,\bs,\bu) = \J(\S_{0,2}(\bu),\bu)
	\leq \liminf_{n\to\infty} \J\,(\S_{\g_n,2}(u^{\g_n}),u^{\g_n})
	\\ & \non
	\leq \liminf_{n\to\infty} \J\,(\S_{\g_n,2}(v),v)
	\leq \lim_{n\to\infty} \J\,(\S_{\g_n,2}(v),v)
	= 	\J\,(\S_{0,2}(v),v)
\end{align}
is satisfied for every fixed $v\in\Uad$, where in the last inequality we 
exploit the continuity of the cost functional with respect to the first
component. \qed

As a consequence of the above proof, we also realize the following corollary.
\Bcor
\label{COR_limJ}
Assume \Propg, \eqref{H1}-\eqref{H7}, \DataApprox, and let $\graffe{\g_n}\subset (0,1]$ 
be a sequence which goes to zero as $n\nearrow\infty$. 
Then, whenever a sequence $\graffe{u^{\g_n}}$ such that $u^{\g_n} \to \bu $ 
weakly star in $L^\infty(Q)$ is given,
we have that 
\Beq
	\non
	\quad \S_{\g_n,2}(u^{\g_n}) \to \S_{0,2}(\bu) \quad \hbox{ weakly star in $\hat{\Y}$}.
\Eeq
Moreover, for every $v\in \UR$, we have that 
\Beq
	\non
	\S_{\g_n,2}(v) \to \S_{0,2}(v) \quad \hbox{strongly in $\L2 V $,}
\Eeq
and, due to the continuity of the cost functional
with respect to the first component, also that
\Beq
	\label{limJ}
	\lim_{n\to\infty} \J(\S_{\g_n,2}(v),v) = \J(\S_{0,2}(v),v) \quad \hbox{ for every $v\in \UR$.}
\Eeq
\Ecor

Below, we prove Theorem \ref{THM_Approximation}, which is the best we can say
as far as the approximation of optimal controls for $(P_0)$ is concerned.
Monotonicity and compactness arguments will be the key arguments to show
this result.

\proof[Proof of Theorem \ref{THM_Approximation}]
Let $\g\in(0,1]$ be fixed, and let $(\bph^\g, \bs^\g, \bu)$ 
be an optimal triple for $\PGtil$ which exists by virtue of Lemma \ref{LM_existencePGtil}.
The boundedness of $\Uad$ and Lemma \ref{LM_regG} yield that
there exist $\ph, \s, u$, and a sequence $\graffe{\g_n}$, which goes to zero as
$n\nearrow\infty$, for which the convergences
\Bsist
	&&
	\bu^{\g_n} \to {u} \ \ \hbox{weakly star in } L^\infty(Q)
	\\  &&
	\bph^{\g_n} \to \ph \ \ \hbox{weakly star in }  \H1 H \cap \L\infty V \cap \L2 W
	\\  &&
	\bs^{\g_n} \to \s \ \ \hbox{weakly star in }  \H1 H \cap \L\infty V \cap \L2 W
\Esist
are satisfied. Moreover, in view of Corollary \ref{COR_limJ},
we also realize that $(\ph,\s,u)$ is an admissible 
triple for $(P_0)$, i.e. that $(\ph,\s)=\S_{0,2}(u)$.

As a matter of fact, monotonicity arguments will allow us to say more. Indeed,
we are able to show that
the limit $u$ coincides with $\bu$, and from the well-posedness of $\S_0$ also that 
$(\ph,\s)=\S_{0,2}(u)=\S_{0,2}(\bu)=(\bph,\bs)$.
As far as $\widetilde{\J}$ is lower semicontinuous and $(\bph,\bs, \bu)$ is optimal for 
$(P_0)$, we deduce that the following inequality holds true
\Bsist
	\label{liminf}
	\non &&
	\liminf_{n\to\infty} \widetilde{\J} (\bph^{\g_n}, \bs^{\g_n}, \bu^{\g_n} )
	\geq \widetilde{\J} (\ph,\s,u) 
	= \J (\ph,\s,u) +\frac 12 \norma{u-\bu}^2_{L^2(Q)}
	\\ && \quad
	\geq \J (\bph,\bs,\bu) +\frac 12 \norma{u-\bu}^2_{L^2(Q)},
\Esist
where we also exploit the definition of the reduced cost functional $\widetilde{\J}$ 
introduced by \eqref{costfunctloc}.
In addition, owing to the optimality of $(\bph^{\g_n}, \bs^{\g_n}, \bu^{\g_n} )$ for 
$(\widetilde{P}_{\g_n})$, we realize that
\Bsist
	\non
	\widetilde{\J} (\bph^{\g_n}, \bs^{\g_n}, \bu^{\g_n} )
	= \widetilde{\J} (\S_{\g_n,2}(\bu^{\g_n} ), \bu^{\g_n} )
	\leq \widetilde{\J} (\S_{\g_n,2}(\bu), \bu )
	\quad \hbox{ for every $n\in\enne$.}
\Esist
Therefore, we pass to the superior limit in both sides
of the above inequality to obtain that
\Bsist
	\label{limsup}
	\limsup_{n\to\infty} \widetilde{\J} (\bph^{\g_n}, \bs^{\g_n}, \bu^{\g_n} )
	\leq
	\widetilde{\J} (\S_{0,2}(\bu), \bu )
	= \widetilde{\J} (\bph,\bs,\bu)
	= \J (\bph,\bs,\bu),
\Esist
where the last equality has been treated invoking \eqref{limJ} and the fact that 
$\widetilde{\J}$ reduces to $\J$ whenever it is considered to act on an 
optimal triple for $(P_0)$. 
Thus, combining inequality \eqref{liminf} and \eqref{limsup} imply that
\Beq
	\frac 12 \norma{u-\bu}^2_{L^2(Q)} = 0,
\Eeq
which consists of the first convergence \eqref{stimaadapteduno}. Moreover, this also establishes
that $\bu=u$ and $(\bph,\bs)=(\ph,\s)$ which prove 
\eqref{stimaadapteddue} and \eqref{stimaadaptedtre} as well.
Finally, due to the above estimates, we also find that
\begin{align}
	\non &
	\J(\bph,\bs,\bu) 
	= \widetilde{\J} (\bph,\bs,\bu)
	= \liminf_{n\to\infty} \widetilde{\J} (\bph^{\g_n}, \bs^{\g_n}, \bu^{\g_n} )
	\\ & \quad \non
	= \limsup_{n\to\infty} \widetilde{\J} (\bph^{\g_n}, \bs^{\g_n}, \bu^{\g_n} )
	= \lim_{n\to\infty} \widetilde{\J} (\bph^{\g_n}, \bs^{\g_n}, \bu^{\g_n} ),
\end{align}
so that the desired convergence \eqref{stimaadaptedquattro} has been shown.
\qed

\section{Optimality results}
\label{SEC_OPTIMALITY}
\setcounter{equation}{0}
This final section is devoted to the check Theorems \ref{THM_stimeadj} 
and \ref{THM_optimality}.
We aim to characterize the optimality conditions
for the initial system \EQ~by passing to the limit in the 
first-order necessary conditions of the approximating problem.
As already underlined, the mathematical analysis that comes out will be quite involved since
a sort of double approximation has to be considered.
To begin with, let us start to investigate the optimality results for the approximating 
system \EQG~which perfectly fits the framework of {\cite{S}}.
As far as numerous problems were already investigated
there, no repetition for the proofs is needed in what follows.

\subsection{The linearized system of the approximating problem}
In the remainder, we require that $\bu\in\UR$ is given, and we denote
$(\bm^\g,\bph^\g,\bs^\g)$ the corresponding solution to \EQG.
At this stage is not so important if $\bu$ is an optimal control or not,
it only matters that $\bu$ is fixed.

For a given $\g\in(0,1]$ and for every $\psi\in L^2(Q)$, the linearized 
system corresponding to \EQG~reads as
follows (c.f. \cite[Sec.~4.2]{S})
\begin{align}
   \a \dt \et^{\g} + \dt \th^{\g} - \Delta \et^{\g} 
  &= P'(\bph^{\g}) (\bs^{\g} - \bm^{\g})\th^{\g} + P(\bph^{\g})(\r^{\g} - \et^{\g})
  \quad \hbox{in $\, Q$}
  \label{EQLinprima}
  \\
  \et^{\g}  &= \b \dt \th^{\g} - \Delta \th^{\g} + g(\g) h''(\bph^{\g}) \th^{\g} + \pi'(\bph^{\g}) \th^{\g}
  \label{EQLinseconda}
  \quad \hbox{in $\,Q$}
  \\
  \dt \r^{\g} - \Delta \r^{\g} & = -P'(\bph^{\g}) (\bs^{\g} - \bm^{\g})\th^{\g} - P(\bph^{\g})(\r^{\g} - \et^{\g}) + \psi
  \label{EQLinterza}
  \quad \hbox{in $\,Q$}
  \\
   \dn \r^{\g} &= \dn \th^{\g} = \dn \et^{\g} = 0
  \quad \hbox{on $\,\Sigma$}
  \label{BCEQLin}
  \\
   \r^{\g}(0) &= \th^{\g}(0) = \et^{\g} (0) = 0
  \quad \hbox{in $\,\Omega$}.
  \label{ICEQLin}
\end{align}
\Accorpa\EQLin EQLinprima ICEQLin

An application of \cite[Thm.~2.4]{S} leads to the result below.
Let us only point out that the symbol $h$ covers a different role there.

\Blem
Assume \Propg, \eqref{H1}-\eqref{H7}, and \DataApprox.
Then, whenever $\psi\in L^2(Q)$ is given, the system \EQLin~admits a unique solution
$(\et^{\g},\th^{\g},\r^{\g})$ that belongs to $\Y$.
\Elem

\subsection{\Frechet~differentiability of $\boldsymbol{\S_\g}$}
Our next goal is concerned with the \Frechet~differentiability of the control-to-state mapping $\S_\g$.
Again, since the system \EQG~was already investigated,
let us recall the obtained result (c.f. \cite[Thm.~2.5]{S}).

\Blem[\bf \Frechet~differentiability of $\S_\g$]
\label{LEM_FreG}
Assume that \Propg, \eqref{H3}-\eqref{H5}, and \DataApprox~are fulfilled. 
Then the control-to-state mapping $\S_\g$, viewed as
a mapping from ${\cal U}_R$ into the state space ${\Y}$, is 
\Frechet~differentiable at $\bar{u}$.
Moreover, for any $\bar{u} \in \UR$ the \Frechet~derivative 
$D\S_\g(\bar{u})$ is a linear and continuous operator 
from $L^2(Q)$ to $\calY$.
Furthermore, for every $ \psi \in L^2(Q)$ we have that
$D\S_\g(\bar{u}) \psi = (\et^\g, \th^\g, \r^\g)$, where $(\et^\g, \th^\g, \r^\g)$ is the unique solution
to \EQLin~associated with $\psi$.
\Elem

\subsection{Optimality conditions for the adapted problem}
In the following, we deal with the optimality conditions for $\PGtil$ that will turn out to be
extremely fruitful in view of the forthcoming asymptotic analysis. At a formal stage, 
we can assert that wherever $\bu^\g$ represents an optimal control 
for $(\widetilde{P}_{\g})$, which exists by
virtue of Lemma \ref{LM_existencePGtil}, the variational inequality,
which characterizes the optimality conditions we are looking for, reads 
\Beq
	\label{abstractcondnecgamma}
	\< D \widetilde{\J}_{\textrm{red},\g}(\bu^{\g}), v - \bu^{\g}> \geq 0 \quad
	\forall v \in \Uad ,
\Eeq
where $D \widetilde{\J}_{\textrm{red},\g}$ denotes the \Frechet~derivative of $\widetilde\J_{\textrm{red},\g}$
and where this latter stands for the reduced cost functional corresponding 
to $\tilde\J$ and it can be defined as made in \eqref{redcost}, while
referred to $\widetilde\J$ instead of $\J$.
Moreover, we can appeal to the chain rule to infer that,
for every fixed $\g\in(0,1]$, we have 
\Beq
	\label{chainrule}
	D \widetilde\J_{\textrm{red},\g}(\bu) = 
	D_{(\bph,\bs)} \widetilde{\J}(\S_{\g,2}(\bu),\bu) \circ D \S_{\g,2}(\bu)
	+ D_{\bu} \widetilde{\J}(\S_{\g,2}(\bu),\bu).
\Eeq
Having Lemma \ref{LEM_FreG} at disposal, a simple application of {\cite[Cor.~2.7]{S}} yields:

\Bcor[\bf First necessary condition]
\label{CORprimanec}
Suppose that the assumptions \Propg, \eqref{H1}-\eqref{H5}, and \DataApprox~are fulfilled.
Let $\g\in(0,1]$ be given, and let $\bu^\g \in \Uad$ be an optimal control for 
$\PGtil$ with its corresponding state
$(\bm^\g,\bph^\g,\bs^\g)$. Then, the necessary condition for the
optimality reads as follows
\begin{align}
  \non
  & \bQ  \intQ (\bph^\g - \phQ)\th^\g 
  + \bO \iO (\bph^\g(T) - \phO)\th^\g(T)
  + \, {b_3}  \intQ (\bs^\g - \sQ)\r^\g 
  \\ & \quad
  + {b_4} \iO (\bs^\g(T) - \sO)\r^\g(T)
  +  \intQ (\bz\bu^\g + (\bu^\g - \bu))(v - \bu^\g)
  \geq 0 \quad \forall v \in \Uad,
  \label{primanec}
\end{align}
where $\th^\g$ and $\r^\g$ are the second and third components of the unique solution 
$(\et^\g,\th^\g,\r^\g)$ to the linearized 
system \EQLin~associated with $ \psi = v - \bu^\g$.
\Ecor
Anyhow, the presence of the linearized variables $\th^\g$ and $\r^\g$ in the above inequality 
is rather unpleasant, thus we try to eliminate them by solving the \socal~adjoint system
that consists of the problem below
\begin{align}
  \non
  & \b \dt q^\g - \dt p^\g + \Delta q^\g - g(\g)h''(\bph^\g)q^\g 
  \\ 
  & \quad 
  - \pi''(\bph^\g) q^\g + P'(\bph^\g)(\bs^\g - \bm^\g)(r^\g-p^\g)= \bQ(\bph^\g - \phQ)
  \quad \hbox{in $\, Q$}
  \label{EQAggprima}
  \\
  & q^\g -\a \dt p^\g - \Delta p^\g + P(\bph^\g)(p^\g -r^\g) = 0
  \label{EQAggseconda}
  \quad \hbox{in $\,Q$}
  \\
  & -\dt r^\g - \Delta r^\g + P(\bph^\g)(r^\g - p^\g) = {b_3} (\bs^\g - \sQ)
  \label{EQAggterza}
  \quad \hbox{in $\,Q$}
  \\
  &\dn q^\g = \dn p^\g = \dn r^\g = 0
  \quad \hbox{on $\Sigma$}
  \label{BCEQAgg}
 \\  
  & p^\g(T) - \b q^\g(T) =\bO(\bph^\g(T) - \phO),  \,\,
  \a p^\g(T) = 0, 
  \,\, r^\g(T) = b_4 (\bs^\g(T) - \sO)
  \quad \hbox{in $\Omega$.}
  \label{ICEQAgg}
\end{align}
\Accorpa\EQAgg EQAggprima ICEQAgg
The well-posedness of the above system has been discussed in \cite{S} where the following
result was proved.
{\Bthm
\label{THMadjoint}
Under the assumptions \Propg, \eqref{H1}-\eqref{H6}, \DataApprox~and
for every fixed $\g\in(0,1]$, the system \EQAgg~has a unique 
solution $(\bar{q}^\g,\bar{p}^\g,\bar{r}^\g)$ that in turn
satisfies the following regularity requirements
\Bsist
	\bar{q}^\g,\bar{p}^\g,\bar{r}^\g \in \H1 H \cap \L\infty V \cap \L2 W \andrea{\subset \C0 V}.
	\label{regadj}
\Esist
\Ethm
This will allow us to obtain the corresponding necessary condition for optimality which reads
as follows.}

\Bthm[\bf Well-posedness and necessary condition]
\label{THM_adjointgamma}
Under the assumptions \Propg, \eqref{H1}-\eqref{H6}, and \DataApprox, 
whenever $\bu^\g$ represents an optimal 
control for $(\widetilde{P}_{\g})$, there holds
\Bsist
  \label{secondanecgamma}
   \intQ (\bar{r}^\g + (\bz \bu^\g + (\bu^\g - \bu))(v - \bu^\g)
  \geq 0 \quad \forall v \in \Uad,
\Esist
where $\bar{r}^\g$ is the unique solution to the adjoint problem introduced by Theorem \ref{THMadjoint}.
\Ethm
At this stage, we would be tempted to pass to the limit 
as $\g \searrow 0$ in the above inequality 
to characterize the necessary conditions for $(P_0)$, but, 
unfortunately, the corresponding mathematical analysis turns out to be more delicate and
the precise description is the purpose of the next paragraph.

\subsection{First-order necessary condition for $\boldsymbol{(P_0)}$}
As sketched in the above lines, we try to recover some optimality conditions
for system \EQ~from inequality \eqref{secondanecgamma} showing that, in a proper sense,
we can pass to the limit as $\g \searrow 0$.
In this direction, some 
compactness properties for the solution to the adjoint problem \EQAgg~need to be shown.
Before moving on, let us introduce a further notation:
in addition to \eqref{QS}, it will be useful to set the backward-in-time
cylinder by setting
\Beq
	\non
	Q_t^T:= \Omega \times [t, T], \quad \hbox{for $0\leq t < T$.}
\Eeq

\proof[Proof of Theorem \ref{THM_stimeadj}]
First, let us show some uniform estimates, with respect to $\g$, that
will allow us to justify the passage to the limit as the parameter goes to zero.

{\bf First estimate:} 
To begin with, let us add to both the members of \eqref{EQAggseconda} the term ${\bar{p}}^\g$. 
Then, let us test \eqref{EQAggprima} by $-\bar{q}^\g$, this new second equation
by $- \dt {\bar{p}}^\g $, and \eqref{EQAggterza} by ${\bar{r}}^\g$. 
Finally, we add these equalities and 
integrate over $Q_t^T$ and by parts to find that
\Bsist
	&& \non
	\frac \b 2 \iO|\bar{q}^\g(t)|^2
	+ \intQtT \dt {\bar{p}}^\g \, \bar{q}^\g
	+ \intQtT |{\nabla \bar{q}^\g}|^2
	+ g(\g) \intQtT h''(\bph^\g)|{{\bar{q}}^\g }| ^2
	- \intQtT \dt {\bar{p}}^\g \, \bar{q}^\g
	\\ && \quad \non
	+ \a \intQtT |{\dt {\bar{p}}^\g}|^2
	+ \frac 12 \iO |\nabla {\bar{p}}^\g(t)|^2
	+ \frac 12 \iO |{\bar{p}}^\g(t)|^2
	+ \frac 12 \iO |{\bar{r}}^\g(t)|^2
	+ \intQtT | {\nabla {\bar{r}}^\g}|^2
	\\ && \quad \non
	=
	\frac \b 2 \iO|\bar{q}^\g(T)|^2
	+ \frac 12 \iO |\nabla {\bar{p}}^\g(T)|^2	
	+ \frac 12 \iO |{\bar{p}}^\g(T)|^2
	+ \frac 12 \iO|{\bar{r}}^\g(T)|^2
	\\ && \quad \non
	- \intQtT \pi''(\bph^\g){\bar{q}^\g \,}^2 
	+ \intQtT P'(\bph^\g)(\bs^\g-\bm^\g)({\bar{r}}^\g-{\bar{p}}^\g)\, \bar{q}^\g
	\\ && \quad \non
	- \intQtT \bQ(\bph^\g-\phQ)\bar{q}^\g
	+ \intQtT P(\bph^\g)({\bar{r}}^\g-{\bar{p}}^\g)\, \dt {\bar{p}}^\g
	- \intQtT {\bar{p}}^\g \, \dt {\bar{p}}^\g
	\\ && \quad \non
	+ \intQtT b_3 (\bs^\g-\sQ){\bar{r}}^\g
	- \intQtT P(\bph^\g)({\bar{r}}^\g-{\bar{p}}^\g){\bar{r}}^\g,
\Esist
where we denote by $I_1,...,I_{11}$ the integrals of the \rhs, in that order.
As regards the \lhs, we point out that 
the second and the fifth integrals cancel out.
Moreover, all the other terms
on that side are nonnegative since we also have that $h''({\bph}^\g)$ is so,
due to the separation result \eqref{sep_g} for $\phG$ and owing to the 
explicit form of $h''$. Therefore, it remains to control the 
\rhs. The assumptions on the final conditions allow us to 
\sfw ly establish that
\Beq
	\non
	\sum_{i=1}^4 |I_i| 
	=
	\frac \bO 2 \iO|\bph^\g(T) - \phO|^2
	+ \frac {b_4}2 \iO|\bs^\g(T) - \sO|^2
	\leq
	c.
\Eeq
Moreover, we have that
\Beq
	\non
	|I_5|+|I_7|+|I_{10}| 
	\leq c 
	+ c  \intQtT | {\bar{q}^\g}|^2 
	+ c  \intQtT | {{\bar{r}}^\g}|^2,
\Eeq
appealing to the Young inequality, to \eqref{H2}, and to \eqref{H5}.
As the remaining integrals are concerned, we compute
\Bsist
	\non &&
	|I_8| + |I_9| + |I_{11}| \leq 
	2\d  \intQtT |   {\dt {\bar{p}}^\g}|^2
	+ \cd \intQtT |  {P(\bph^\g)({\bar{r}}^\g-{\bar{p}}^\g)}|^2
	+ \cd \intQtT |  {{\bar{p}}^\g}|^2
	+c \intQtT |  {{\bar{r}}^\g}|^2
	\\ \non && \quad
	\leq 
	2\d \intQtT |  {\dt {\bar{p}}^\g}|^2
	+ \cd \intQtT |  {{\bar{p}}^\g}|^2
	+ \cd \intQtT |  {{\bar{r}}^\g }|^2,
\Esist
where we have applied \eqref{H4} and \eqref{young}.
Moreover, we obtain from \eqref{H4}, \eqref{H5}, \eqref{reg_gamma_ms},
\eqref{sep_g}, the Sobolev embedding \eqref{sobolev_emb}, Young's and \Holder's inequality that
\Bsist
	\non
	&& |I_6| \leq
	c  \intQtT |   {P'(\bph^\g)(\bs^\g-\bm^\g) \bar{q}^\g}|^2
	+ c  \intQtT |   {{\bar{r}}^\g-{\bar{p}}^\g}|^2
	\\ \non && \quad \leq 
	c \intQtT ( | {\bs^\g }|^2 + | {\bm^\g}| ^2 )|\bar{q}^\g | |\bar{q}^\g |
	+ c  \intQtT |  {\bar{r}^\g }|^2
	+ c  \intQtT |  {\bar{p}^\g}|^2	
	\\ \non && \quad
	\leq 
	c \inttt (\norma{\bs^\g}_6^2 + \norma{\bm^\g}_6^2)\norma{\bar{q}^\g}_6 \norma{\bar{q}^\g}_2
	+ c \intQtT | {\bar{r}^\g}|^2
	+ c \intQtT | {\bar{p}^\g}	|^2
	\\ \non && \quad
	\leq 
	c \inttt (\norma{\bs^\g}_V^2 + \norma{\bm^\g}_V^2)\norma{\bar{q}^\g}_V \norma{\bar{q}^\g}_H
	+ c \intQtT | {\bar{r}^\g}|^2
	+ c \intQtT | {\bar{p}^\g}	|^2
	\\ \non && \quad
	\leq 
	\frac 12 \intQtT (|\bar{q}^\g|^ 2+ |\nabla \bar{q}^\g|^2)
	+ c \inttt (\normaV{\bs^\g}^4 + \normaV{\bm^\g}^4)\norma{\bar{q}^\g}_H^2
	+ c \intQtT | {\bar{r}^\g}|^2
	+ c \intQtT | {\bar{p}^\g}	|^2
	\\ \non && \quad
	\leq
	\frac 12 \intQtT | {\nabla \bar{q}^\g}|^2
	+c \intQtT |  {\bar{r}^\g}|^2
	+ c \intQtT | {\bar{p}^\g}|^2
	+ c \intQtT | {\bar{q}^\g}|^2.
\Esist
Then, upon collecting the above estimates, 
we fix $0<\d <\a/2$ and apply the backward-in-time Gronwall lemma to infer that
\Bsist
	\non &&
	\norma {\bar{q}^\g}_{\L\infty H \cap \L2 V}
	+ \norma {{\bar{p}}^\g}_{\H1 H \cap \L\infty V}
	+ \norma {{\bar{r}}^\g}_{\L\infty H \cap \L2 V }
	\\ && \quad
	+ g(\g) \intQ h''(\bph^\g)|{{\bar{q}}^\g }| ^2
	\leq c.
	\label{stima_uno}
\Esist

{\bf Second estimate:}
We proceed multiplying \eqref{EQAggseconda} by $\Delta {\bar{p}}^\g$.
Using the Young inequality and the above estimate, we deduce that
\Beq
	\non
	\norma{\Delta {\bar{p}}^\g}_{\L2 H} \leq c,
\Eeq
and, accounting for the elliptic regularity theory
for homogeneous Neumann boundary problems, also that
\Beq
	\label{stima_due}
	\norma{{\bar{p}}^\g}_{\L2 W} \leq c.
\Eeq

{\bf Third estimate:}
Next, let us rewrite the equation \eqref{EQAggterza} as follows:
\Beq
  -\dt \bar{r}^\g - \Delta \bar{r}^\g = {b_3} (\bs^\g - \sQ) 
  - P(\bph^\g)(\bar{r}^\g - \bar{p}^\g)=:f.
\Eeq
Owing to the above estimates, \eqref{H4}, and using the fact that $\bph^\g$, 
as solution to \EQ, satisfies \eqref{reg_gamma_ph},
it is easy to realize that the forcing term $f \in\L2 H$.
Hence, since it reads as a backward-in-time parabolic equation, on account
for the boundary condition \eqref{BCEQAgg} and for 
the regularity of the final datum \eqref{ICEQAgg}, we obtain that
\Beq
	\label{stima_tre}
	\norma{{\bar{r}}^\g}_{\H1 H \cap \L\infty V \cap \L2 W} \leq c.
\Eeq

{\bf Fourth estimate:}
Now, let us point out that whenever $v\in {{\cal W}_0(0,T)}$ is given, a simple
application of the integration by parts formula and the last of \eqref{ICEQAgg} lead to
\Beq
	\non
	\b\,_{{\cal W}_0^*} \< \dt \bar{q}^\g, v >_{{\cal W}_0} 
	= 
	- \b \ioT \< \dt v, \bar{q}^\g >
	- \bO \iO (\bph^\g(T) - \phO)  v(T),
\Eeq	
which, in turn, implies that
\Beq
	\non
	\Big|\,_{{\cal W}_0^*} \< \dt \bar{q}^\g, v >_{{\cal W}_0} \Big|
	\leq 
	\b  \norma{\dt v}_{\L2 \Vp} \, \norma{\bar{q}^\g}_{\L2 V}
	+ \norma{\bO (\bph^\g(T) - \phO)}_H \, \norma{v(T)}_H
	\leq c \norma{v}_{{\cal W}_0},
\Eeq	
where the previous estimates and the continuous embedding
${\cal W}_0 (0,T) \subset \C0 H $ are taken into account.
Dividing both sides by $\norma{v}_{{\cal W}_0}$ and passing to the supremum, 
we conclude that there exists a positive constant $c$ such that
\Beq
	\norma{\dt \bar{q}^\g}_{{\cal W}_0^*}
	\leq c.
	\label{stima_quattro}
\Eeq

{\bf Fifth estimate:}
We are left with the task to control $\lam^\g$ which was introduced by
\eqref{deflambda}. First, let us rewrite the equation \eqref{EQAggprima} as
\Bsist
	\non
	&& \lam^\g= g(\g)h''(\bph^\g)\bar{q}^\g 
	= 
	\b \dt \bar{q}^\g 
	- \dt \bar{p}^\g 
	+ \Delta \bar{q}^\g 
	- \pi''(\bph^\g) \bar{q}^\g 
	\\ && \quad \non
	+ P'(\bph^\g)(\bs^\g - \bm^\g)(\bar{r}^\g-\bar{p}^\g) - \bQ(\bph^\g - \phQ).
\Esist
Then, we consider the duality pairing between $ \lam^\g$ and an arbitrary function 
$v \in {\cal W}_0(0,T)$ to get
\begin{align}
	\non &
	\,_{{\cal W}_0^*} \< \lam^\g, v >_{{\cal W}_0} = 
	 \b \,_{{\cal W}_0^*}\<   \dt \bar{q}^\g, v>_{{\cal W}_0^*}
	- \intQ \dt \bar{p}^\g  v 
	- \intQ \nabla \bar{q}^\g \cdot \nabla v
	- \intQ \pi''(\bph^\g) \bar{q}^\g v
	\\ & \quad 
	+ 	\intQ P'(\bph^\g)(\bs^\g - \bm^\g)(\bar{r}^\g-\bar{p}^\g) v
	-  \intQ\bQ(\bph^\g - \phQ) v,
	\label{lam_conf}
\end{align}	
where also \eqref{BCEQAgg} and \eqref{ICEQAgg} are taken into account.
Hence, due to \eqref{H2}, \eqref{H5} and to the previous estimates,
we claim that the above inequality implies that
\Bsist
	&&
	\label{lam_wo}
	\Big| {}_{{\cal W}_0^*} \< \lam^\g, v >_{{\cal W}_0}\Big| 
	\leq
	c \norma{ v }_{{\cal W}_0},
\Esist
for a positive constant $c$.
Indeed, bearing in mind \accorpa{stima_uno}{stima_due}, 
and \accorpa{stima_tre}{stima_quattro}, it is clear that
the only term which deserves further investigation is the fifth product on the \rhs~of
\eqref{lam_conf} that we denote by $\widetilde{I}$.
Then, let us separately prove how it can be controlled. 
{By virtue of the \Holder~inequality, the boundedness of $P'$ and the Sobolev embedding
\eqref{sobolev_emb}, we conclude that
\begin{align}
	\non &
	|\widetilde{I}|
	\leq
	c \intQ |\bs^\g - \bm^\g| |\bar{r}^\g-\bar{p}^\g| |v|
	\leq
	c \ioT \norma{\bs^\g - \bm^\g}_4 \norma{\bar{r}^\g-\bar{p}^\g}_2 \norma{v}_4
	\\ & \quad
	\leq
	c \ioT \norma{\bs^\g - \bm^\g}_V \norma{\bar{r}^\g-\bar{p}^\g}_H \norma{v}_V
	\leq
	c \norma{\bar{r}^\g-\bar{p}^\g}_{\L2 H} \norma{v}_{\L2 V},
\end{align}
where in the last inequality we also appeal to the fact that $\bs^\g$ and $\bm^\g$,
as solutions to \EQG, satisfy \eqref{reg_gamma_ph} and \eqref{reg_gamma_ms} and therefore they
both belong to $\L\infty V$.
So that \eqref{lam_wo} is shown.}
Lastly, dividing both sides of \eqref{lam_wo} by $\norma{v}_{{\cal W}_0}$ and passing to the supremum,
we obtain that there exists a positive constant $c$ such that
\Beq
	\label{stima_cinque}
	\norma{\lam^\g}_{{\cal W}_0^*}
	\leq c.
\Eeq
To conclude, the application of the aforementioned estimates \accorpa{stima_uno}{stima_due}, 
\accorpa{stima_tre}{stima_quattro}, and \eqref{stima_cinque}, lead to infer
the existence of some $q,p,r$ and $\lam$ such that, as $n \nearrow \infty$, 
\Conv, are verified.
\qed

We will see that the above result will be sufficient to pass to the limit as $\g \searrow 0$, at least 
in a suitable weak framework, as rigorously described in Theorem \ref{THM_optimality}.

\proof[Proof of Theorem \ref{THM_optimality}]
Owing to the uniform estimates pointed out in Theorem \ref{THM_stimeadj}, 
we realize that there exists some subsequence, which is again indexed by $n$, for which we 
can pass to the limit in the inequality \eqref{secondanecgamma}
to obtain that, for some limit $ r $, \eqref{nec_final} is satisfied.

Thus, we are led to show
that such a limit solves a suitable adjoint problem for \EQ, and
eventually prove some additional features on this system.
To do that, let us multiply the equations \eqref{EQAggprima} by 
an arbitrary $v\in{\cal W}_0(0,T)$ and integrate over $Q$.
Accounting for the boundary conditions \eqref{BCEQAgg}, and for the final ones \eqref{ICEQAgg},
the application of the integration by parts, 
leads to the following problem consisting of a variational formulation
\Bsist
	\non
	& \ioT \<  \dt v , \bar{p}^{\g_n} - \b \bar{q}^{\g_n}>	
	- {}_{{\cal W}_0^*}\< \lam^{\g_n} , v >_{{\cal W}_0}	
	+ \intQ \nabla \bar{q}^{\g_n} \cdot \nabla v
	\\ & \quad  \non
	- \intQ \pi''(\bph^{\g_n}) \bar{q}^{\g_n} v
	+ \intQ P'(\bph^{\g_n})(\bs^{\g_n} - \bm^{\g_n})({\bar{r}}^{\g_n}-{\bar{p}}^{\g_n})\, v
	\\ & \quad \label{varfor}
	= 
	\intQ \bQ (\bph^{\g_n} - \phQ) v
	+ \iO \bO(\bph^{\g_n}(T) - \phO) v(T)
	\quad \hbox{for every $v \in {{\cal W}_0(0,T)},$}
\Esist
combined with the system
\begin{align}
	{\bar{q}}^{\g_n} -\a \dt {\bar{p}}^{\g_n} - \Delta {\bar{p}}^{\g_n} 
	+ P(\bph^\g)({\bar{p}}^{\g_n} -{\bar{r}}^{\g_n})& =   0
    \quad \hbox{in $\,Q$}
    \\
    -\dt {\bar{r}}^{\g_n} 
    - \Delta {\bar{r}}^{\g_n} + P(\bph^\g)({\bar{r}}^{\g_n} - {\bar{p}}^{\g_n})  &= {b_3} (\bs^\g - \sQ)
    \quad \hbox{in $\,Q$}
    \\
    \dn {\bar{p}}^{\g_n} &= \dn {\bar{r}}^{\g_n} = 0
    \quad \hbox{on $\Sigma$}
	\\  \quad
	\a {\bar{p}}^{\g_n}(T) &= 0, 
	\,\, {\bar{r}}^{\g_n}(T) = b_4 (\bs^\g(T) - \sO)
	\quad \hbox{in $\Omega$.}
	\label{finale_adj}
\end{align}
With the convergences \Conv~at our disposal, we would like to let $n \nearrow \infty $ to 
show that in the limit we get \ADJfin.
Let us recall that, as pointed out by \accorpa{strongconvuno}{strongconvdue}, we have that
\Bsist
	\non
	&&\bph^{\g_n} \to \bph , \quad \bs^{\g_n} \to \bs, \quad \bm^{\g_n} \to \bm
	\quad \hbox{strongly in }  \C0 H \cap \L2 V. 
\Esist
In a similar fashion, we immediately deduce from \Conv~that, up to 
not relabeled subsequence, there holds 
\Bsist
	\non
	&&\bar{p}^{\g_n} \to p \ \ \hbox{strongly in }  \C0 H \cap \L2 V 
	\\  && \non
	\bar{r}^{\g_n} \to r \ \ \hbox{strongly in }  \C0 H \cap \L2 V.
\Esist
Therefore, even though some terms in \accorpa{varfor}{finale_adj} possess a strongly nonlinear
behavior, recalling that $\pi''$, $P$ and $P'$ are continuous, we have that
the above strong convergences allow us to let $n\nearrow \infty$ to obtain that in the limit
\ADJfin~is satisfied.
This is the sense
in which we can state that the limit $(q,p,r)$ and $\lam$ enjoy an adjoint system
corresponding to \EQ. 

To conclude, let us provide some additional properties which characterize
such a limit.
First, let us multiply $ \lam^{\g_n}$, which is defined by 
\eqref{deflambda}, by $\bar{q}^{\g_n} $ and 
integrate over $Q$ to realize that
\Bsist
	\non
	 \intQ \lam^{\g_n} \bar{q}^{\g_n}
	=  \intQ g(\g_n)h''(\bph^{\g_n})|\bar{q}^{\g_n}|^2\geq 0
	\quad  \hbox{for every $n \in \enne$},
\Esist
where we used that $g$ is nonnegative,
that $h''(s)=\frac 2{1-s^2}$, and the fact that ${\bph}^{\g_n}$, as a solution to 
\EQG, verifies the separation result \eqref{sep_g} form which we infer that $ h''(\bph^{\g_n}) > 0$. 
Then, by passing to the inferior limit, as $n\nearrow\infty$, we 
conclude the first condition \eqref{slackness_uno}.
Lastly, let us show a limit behavior that should suggest the fact that,
in the limit, $\lam^{\g_n}$ tends to concentrate in the sets where $|{\bph}| = 1$.
We now multiply $ \lam^{\g_n}$ by $(1-(\bph^{\g_n})^2) \Phi$, for an arbitrary regular function
$\Phi$ which vanishes at zero, integrate over $Q$,
and then pass to the limit as $n\nearrow\infty$ to find 
\Bsist
	\non
	L_{\Phi}:= \lim_{n\to\infty}  \intQ \lam^{\g_n}(1-(\bph^{\g_n})^2) \, \Phi.
\Esist
Moreover, exploiting the definition of $\lam^{\g_n}$, the explicit expression
of $h''$, and accounting for the asymptotic properties \Propg~we are assuming, we realize that 
\Beq
	\non
	L_{\Phi} = 2\, \lim_{n\to\infty}  \intQ g(\g_n) \, \bar{q}^{\g_n} \, \Phi=0,
\Eeq
which is exactly the condition \eqref{slackness_due} we are going to prove.
\qed

\section{Appendix}
\label{APPENDIX}
\setcounter{equation}{0}

Here, we focus on showing a possible way to construct an approximating family of data
which fulfills requirements \DataApprox. 
\Blem
Let $(\m_0,\ph_0,\s_0)$ be a triplet belonging to $V \times V \times V$.
Then, there exists an approximating family $\graffe{(\mG_0,\phG_0,\sG_0)}_\g$ which satisfies
all the following properties 
\Bsist
	& \non
	(\mG_0,\phG_0,\sG_0) \in 
		(\Huno \cap L^\infty(\Omega)) \times W \times \Huno 
		\quad \forall \g\in(0,1],
	\\ \non & 
	|\phG_0| \leq 1- \g/2 \quad \aeO   \quad \forall \g\in(0,1],
	\\	& \non
	(\mG_0,\phG_0,\sG_0)  \to (\m_0,\ph_0,\s_0) \quad \hbox{strongly in } V \times V \times V 
			\quad \hbox{as } \g\searrow 0.
\Esist
\Elem
\proof
Obviously, since only existence is stated above, we only show one possible way to proceed.
As the first and third variables are concerned, the choices are quite natural. 
As $\s_0^\g$ we can \sfw ly take, for every $\g$, $\s_0$ itself.
Secondly, it is natural to choose as $\m_0^\g$ a suitable truncation of 
$\m_0$, since we want $\m_0^\g$ to remain bounded in $V$ and to be uniformly bounded as well. 
So, we take $\m_0^\g$ as the truncation at level $1/\g$ of $\m_0$,
namely the function defined by
\Beq
	\non
	\m_0^\g := \begin{cases} 1/\g &\hbox{ if } \m_0  > 1/\g
				\\ \m_0^\g  &\hbox{ if } |\m_0^\g| \leq 1/\g.
				\\ -1/\g  &\hbox{ if } \m_0  < - 1/\g
    \end{cases}
\Eeq
It is now a standard argument to check that $\m_0^\g \to \m_0$ strongly in $V$.
To conclude, let us face the remaining term which will require more attention. 
To recover the zero normal derivative condition, we are tempted to choose as $\ph_0^\g$ the solution to
the following homogeneous Neumann boundary problem
\Beq
 \non
 \begin{cases} 
 	\ph_0^\g - \g \Delta \ph_0^\g &= \ph_0 \quad \hbox{in $\Omega$}
 	\\
 	\dn \ph_0^\g &= 0 \quad \hbox{on $\Gamma$},
 \end{cases}
\Eeq
in which $\ph_0$ appears as forcing term.
Owing to the regularity of $\ph_0$, for every $\g \in (0,1]$, it follows from standard results
that there exists a unique solution $\ph_0^\g $ which belongs to $W$.
Let us mention that in this framework one can also verify 
that $\ph_0^\g \to \ph_0$ strongly in $V$.
However, we also care to control the absolute value of 
$\ph_0^\g$ by $1-\g/2$, which is not guaranteed if we proceed this way. 

Thus, we will follow a similar technique, but first we need to truncate the data 
in order to properly control the absolute value of the approximating variable.
For convenience, let us denote as $\tilde{\ph}_0^\g$ the truncation of $\ph_0$ at level
$1-\g/2$. Then, we consider the following problem which strictly resembles the above one
\Beq
 \begin{cases} 
 	\ph_0^\g - \g \Delta \ph_0^\g & = \tilde{\ph}_0^\g \quad \hbox{in $\Omega$}
 	\\
 	\dn \ph_0^\g & = 0\quad \hbox{on $\Gamma$}.
 	\label{pertsing}
 \end{cases}
\Eeq
Similarly, it is easy deduce that, for every $\g\in (0,1]$, the unique solution $\ph_0^\g$
to \eqref{pertsing} belongs to $W$.
Moreover, by adding to both sides of the first equation of \eqref{pertsing} the term $1-\g/2$, we 
arrive at the identity
\Beq
	\non
	(\ph_0^\g +1-\g/2) - \g \Delta (\ph_0^\g \g+1-\g/2)= \tilde{\ph}_0^\g +1-\g/2,
\Eeq
where the \rhs~is positive by construction of $\tilde{\ph}_0^\g$. 
Hence, the maximum principle yields that 
\Beq
	\non
	\ph_0^\g +1-\g/2 \geq 0,
	\hbox{ which implies }
	\ph_0^\g \geq -1+ \g/2. 
\Eeq
By repeating the same strategy adding to both sides the term $-1+\g/2$,
it is easy to conclude that $|\ph_0^\g| \leq 1-\g/2$.
Finally, we are left with the task of showing that $\ph_0^\g  \to \ph_0$ strongly in $V$.
In this direction, we multiply the first equation of 
\eqref{pertsing} by $\ph_0^\g$ and integrate over $\Omega$.
Using the boundary condition and estimating the \rhs~by mean of the Young inequality, 
we discover that
\Beq
	\non
	\iO |\ph_0^\g|^2
	+ \g \iO |\nabla \ph_0^\g|^2
	\leq 
	\frac 12 \iO |\ph_0^\g|^2
	+ \frac 12 \iO |\tilde{\ph}_0^\g|^2,
\Eeq
and rearranging the terms, we obtain the boundedness of $\graffe{\ph_0^\g}_\g$ in $V$.
Thus, from weak compactness arguments, we immediately realize that 
$\ph_0^\g  \to \ph_0$ at least weakly in $V$.
Now, we multiply the first equation by $- \Delta \ph_0^\g$ and integrate by parts to get
\Beq
	\non
	\iO |\nabla \ph_0^\g|^2
	+ \g \iO |\Delta \ph_0^\g|^2
	= 
	\iO \nabla \tilde{\ph}_0^\g \cdot \nabla \ph_0^\g,
\Eeq
where we also account for the homogeneous boundary condition for $\ph_0^\g$.
Due to the Young inequality
we control the \rhs~as follows
\Beq
	\non
	\iO |\nabla \ph_0^\g|^2
	+ \g \iO |\Delta \ph_0^\g|^2
	\leq 
	\frac 12 \iO |\nabla \ph_0^\g|^2
	+ \frac 12 \iO |\nabla \tilde{\ph}_0^\g|^2,
\Eeq
and, adding this latter with the above estimate and rearranging the terms, we obtain that
\Beq
	\non
	\frac 12 \iO |\ph_0^\g|^2
	+ \frac 12 \iO |\nabla \ph_0^\g|^2
	+ \g \iO |\nabla \ph_0^\g|^2
	+ \g \iO |\Delta \ph_0^\g|^2
	\leq 
	 \frac 12 \iO |\tilde{\ph}_0^\g|^2
	+ \frac 12 \iO |\nabla \tilde{\ph}_0^\g|^2,
\Eeq
from which we realize that 
$\norma{\ph_0^\g}^2_V \leq \norma{\tilde{\ph}_0^\g}^2_V$,
and since that inequality continue to hold if we pass to the superior limit, 
we realize that actually $\ph_0^\g  \to \ph_0$ strongly in $V$, as $\g\searrow 0$,
which is the desired conclusion.
\qed 


\vspace{3truemm}

\Begin{thebibliography}{10} 
\footnotesize

\bibitem{Agosti}
A. Agosti, P.F. Antonietti, P. Ciarletta, M. Grasselli and M. Verani,
A Cahn--Hilliard--type equation with application to tumor growth dynamics, 
{\it Math. Methods Appl. Sci.}, {\bf 40} (2017), 7598-–7626.

\bibitem{BARBU}
V. Barbu,
Necessary conditions for nonconvex distributed control problems governed by 
elliptic variational inequalities, 
{\it J. Math. Anal. Appl.} {\bf 80} (1981), 566-597.

\bibitem{BRZ}
H. Brezis,
``Op\'erateurs maximaux monotones et semi-groupes de contractions dans les
espaces de Hilbert'', North-Holland Math. Stud. {\bf 5}, North-Holland, Amsterdam, (1973).

\bibitem{CRW}
C. Cavaterra, E. Rocca and H. Wu,
Long--time Dynamics and Optimal Control of a Diffuse Interface Model for Tumor Growth,
{\it Appl. Math. Optim.} (2019), https://doi.org/10.1007/s00245-019-09562-5.

\bibitem{CFGS}
P. Colli, M.H. Farshbaf-Shaker, G. Gilardi and J. Sprekels,
Optimal boundary control of a viscous Cahn–-Hilliard system
with dynamic boundary condition and double obstacle potentials,
{\it SIAM J. Control Optim.} {\bf 53} (2015), 2696--2721.

\bibitem{CFS}
P. Colli, M.H. Farshbaf-Shaker and J. Sprekels,
A deep quench approach to the optimal control of an Allen–-Cahn equation
with dynamic boundary conditions and double obstacles,
{\it Appl. Math. Optim.} {\bf 71} (2015), 1-24.

\bibitem{CGH}
P. Colli, G. Gilardi and D. Hilhorst,
On a Cahn-Hilliard type phase field system related to tumor growth,
{\it Discrete Contin. Dyn. Syst.} {\bf 35} (2015), 2423-2442.

%
\bibitem{CGRS_VAN}
P. Colli, G. Gilardi, E. Rocca and J. Sprekels,
Vanishing viscosities and error estimate for a Cahn–-Hilliard type phase field system related to tumor growth,
{\it Nonlinear Anal. Real World Appl.} {\bf 26} (2015), 93-108.

\bibitem{CGRS_OPT}
P. Colli, G. Gilardi, E. Rocca and J. Sprekels,
Optimal distributed control of a diffuse interface model of tumor growth,
{\it Nonlinearity} {\bf 30} (2017), 2518-2546.

\bibitem{CGRS_ASY}
P. Colli, G. Gilardi, E. Rocca and J. Sprekels,
Asymptotic analyses and error estimates for a \CH~type phase field system modeling tumor growth,
{\it Discrete Contin. Dyn. Syst. Ser. S} {\bf 10} (2017), 37-54.

%
%

\bibitem{CGS_DQ}
P. Colli, G. Gilardi and J. Sprekels,
Optimal velocity control of a convective Cahn-–Hilliard system with double obstacles
and dynamic boundary conditions: a `deep quench' approach.
{\it J. Convex Anal.} {\bf 26} (2019).

\bibitem{CGS_nonlocal}
P. Colli, G. Gilardi and J. Sprekels,
Distributed optimal control of a nonstandard nonlocal phase field system with double obstacle potential,
{\it Evol. Equ. Control Theory} {\bf 6} (2017), 35--58.

%
%
%

\bibitem{CS_doubleobstacle}
P. Colli and J. Sprekels,
Optimal boundary control of a nonstandard Cahn-–Hilliard system
with dynamic boundary condition and double obstacle inclusions,
in {\sl Solvability, Regularity, Optimal Control of Boundary Value Problems for PDEs},
P.~Colli, A.~Favini, E.~Rocca, G.~Schimperna, J.~Sprekels~(ed.),
{\it Springer INdAM Series} {\bf 22}, Springer, Milan, 2017, 151--182.

\bibitem{CLLW}
V. Cristini, X. Li, J.S. Lowengrub, S.M. Wise,
Nonlinear simulations of solid tumor growth using a mixture model: invasion and branching.
{\it J. Math. Biol.} {\bf 58} (2009), 723–-763.

\bibitem{CL}
V. Cristini, J. Lowengrub,
Multiscale Modeling of Cancer: An Integrated Experimental and Mathematical
{\it Modeling Approach. Cambridge University Press}, Leiden (2010).

\bibitem{DFRGM}
M. Dai, E. Feireisl, E. Rocca, G. Schimperna, M. Schonbek,
Analysis of a diffuse interface model of multispecies tumor growth,
{\it Nonlinearity\/} {\bf  30} (2017), 1639--1658.

\bibitem{EK_ADV}
M. Ebenbeck and P. Knopf,
Optimal control theory and advanced optimality conditions for a diffuse interface model of tumor growth 
{\it preprint arXiv:1903.00333 [math.OC],} (2019), 1--34.

\bibitem{EK}
M. Ebenbeck and P. Knopf,
Optimal medication for tumors modeled by a Cahn--Hilliard--Brinkman equation,
{\it preprint arXiv:1811.07783 [math.AP],} (2018), 1--26.

\bibitem{EGAR}
M. Ebenbeck and H. Garcke,
Analysis of a Cahn–-Hilliard–-Brinkman model for tumour growth with chemotaxis.
{\it J. Differential Equations,} (2018) https://doi.org/10.1016/j.jde.2018.10.045.

\bibitem{FGR}
S. Frigeri, M. Grasselli, E. Rocca,
On a diffuse interface model of tumor growth,
{\it  European J. Appl. Math.\/} {\bf 26 } (2015), 215-243. 

%
\bibitem{FLRS}
S. Frigeri, K.F. Lam, E. Rocca, G. Schimperna,
On a multi-species Cahn--Hilliard-Darcy tumor growth model with singular potentials,
{\it Comm. Math. Sci.\/} {\bf 16 (3)}, (2018), 821-–856. 

\bibitem{FRL}
S. Frigeri, K.F. Lam and E. Rocca,
On a diffuse interface model for tumour growth with non-local interactions and degenerate 
mobilities,
In {\sl  Solvability, Regularity, and Optimal Control of Boundary Value Problems for PDEs},
P. Colli, A. Favini, E. Rocca, G. Schimperna, J. Sprekels (ed.),
{\it Springer INdAM Series,} {\bf 22}, Springer, Cham, 2017.

\bibitem{GARL_1}
H. Garcke and K. F. Lam,
Well-posedness of a Cahn-–Hilliard–-Darcy system modelling tumour
growth with chemotaxis and active transport,
{\it European. J. Appl. Math.} {\bf 28 (2)} (2017), 284-–316.
\bibitem{GARL_2}
H. Garcke and K. F. Lam,
Analysis of a Cahn--Hilliard system with non--zero Dirichlet 
conditions modeling tumor growth with chemotaxis,
{\it Discrete Contin. Dyn. Syst.} {\bf 37 (8)} (2017), 4277-4308.
\bibitem{GARL_3}
H. Garcke and K. F. Lam,
Global weak solutions and asymptotic limits of a Cahn--Hilliard--Darcy system modelling tumour growth,
{\it AIMS Mathematics} {\bf 1 (3)} (2016), 318-360.
\bibitem{GARL_4}
H. Garcke and K. F. Lam,
On a Cahn--Hilliard--Darcy system for tumour growth with solution dependent source terms, 
in {\sl Trends on Applications of Mathematics to Mechanics}, 
E.~Rocca, U.~Stefanelli, L.~Truskinovski, A.~Visintin~(ed.), 
{\it Springer INdAM Series} {\bf 27}, Springer, Cham, 2018, 243-264.
\bibitem{GARLR}
H. Garcke, K. F. Lam and E. Rocca,
Optimal control of treatment time in a diffuse interface model of tumor growth,
{\it Appl. Math. Optim.} {\bf 78}(3) (2018), {495--544}.
\bibitem{GAR}
H. Garcke, K. F. Lam, R. N\"urnberg and E. Sitka,
A multiphase Cahn--Hilliard--Darcy model for tumour growth with necrosis,
{\it Math. Models Methods Appl. Sci.}  {\bf 28 (3)} (2018), 525-577.


\bibitem{HDPZO}
A. Hawkins-Daarud, S. Prudhomme, K.G. van der Zee, J.T. Oden,
Bayesian calibration, validation, and uncertainty quantification of diffuse 
interface models of tumor growth. 
{\it J. Math. Biol.} {\bf 67} (2013), 1457–-1485. 

\bibitem{HDZO}
A. Hawkins-Daruud, K. G. van der Zee and J. T. Oden, Numerical simulation of
a thermodynamically consistent four-species tumor growth model, 
{\it Int. J. Numer. Math. Biomed. Engng.} {\bf 28} (2011), 3–-24.

\bibitem{HKNZ}
D. Hilhorst, J. Kampmann, T. N. Nguyen and K. G. van der Zee, Formal asymptotic
limit of a diffuse-interface tumor-growth model, 
{\it Math. Models Methods Appl. Sci.} {\bf 25} (2015), 1011-–1043.


%
%

\bibitem{Lions_OPT}
J. L. Lions,
Contr\^ole optimal de syst\`emes gouverne\'s par des equations aux d\'eriv\'ees partielles,
Dunod, Paris, 1968.

\bibitem{Mir_CH}
A. Miranville,
The Cahn-Hilliard equation and some of its variants, 
{\it AIMS Mathematics,} {\bf 2} (2017), 479–-544.

\bibitem{S}
A. Signori,
Optimal distributed control of an extended model of tumor 
growth with logarithmic potential. 
{\it Appl. Math. Optim.} (2018), https://doi.org/10.1007/s00245-018-9538-1.

\bibitem{Simon}
J. Simon,
{Compact sets in the space $L^p(0,T; B)$},
{\it Ann. Mat. Pura Appl.~(4)\/} 
{\bf 146} (1987) 65--96.


\bibitem{SW}
J. Sprekels and H. Wu,
Optimal Distributed Control of a Cahn-–Hilliard–-Darcy System with Mass Sources,
{\it Appl. Math. Optim.} (2019), https://doi.org/10.1007/s00245-019-09555-4.

\bibitem{Trol}
F. Tr\"oltzsch,
Optimal Control of Partial Differential Equations. Theory, Methods and Applications,
{\it Grad. Stud. in Math.,} Vol. {\bf 112}, AMS, Providence, RI, 2010.

\bibitem{WZZ}
X. Wu, G.J. van Zwieten and K.G. van der Zee, Stabilized second-order splitting
schemes for \CH~models with applications to diffuse-interface tumor-growth models, 
{\it Int. J. Numer. Meth. Biomed. Engng.} {\bf 30} (2014), 180-203.

\End{thebibliography}

\End{document}

\bye